\documentclass[oneside, 12pt]{amsart}

\usepackage{amscd, amssymb, amsmath, mathrsfs}
\usepackage[english]{babel}
\usepackage{booktabs}
\usepackage[all,cmtip]{xy}
\usepackage[colorlinks, linktocpage, citecolor = blue, linkcolor = blue]{hyperref}

\newcommand\mylabel[1]{\label{#1}}

\newtheorem{theorem}{Theorem}[section]
\newtheorem*{maintheorem}{Theorem}
\newtheorem{lemma}[theorem]{Lemma}
\newtheorem{proposition}[theorem]{Proposition}
\newtheorem{corollary}[theorem]{Corollary}

\theoremstyle{definition}
\newtheorem{definition}[theorem]{Definition}

\newtheorem*{acknowledgement}{Acknowledgement}

\theoremstyle{remark}

\newcommand{\ZZ}	{\mathbb{Z}}

\newcommand{\CC}	{\mathbb{C}}

\newcommand  {\shC}     {\mathscr{C}}

\newcommand  {\shF}     {\mathscr{F}}
\newcommand  {\shG}     {\mathscr{G}}

\newcommand  {\shI}     {\mathscr{I}}

\newcommand  {\shL}     {\mathscr{L}}
\newcommand  {\shP}     {\mathscr{P}}
\newcommand  {\shR}     {\mathscr{R}}
\newcommand  {\shS}     {\mathscr{S}}
\newcommand  {\shT}     {\mathscr{T}}

\newcommand  {\foC}     {\mathfrak{C}}

\newcommand  {\foG}     {\mathfrak{G}}
\newcommand  {\foH}     {\mathfrak{H}}

\newcommand  {\foR}     {\mathfrak{R}}

\newcommand  {\foT}     {\mathfrak{T}}

\newcommand	{\Ab}	{\text{\rm Ab}}

\newcommand	{\alt}	{\text{\rm alt}}

\newcommand  {\can}     { \text{\rm can}}

\newcommand  {\cH}      {\check{H}}

\newcommand  {\Ch}      {\operatorname{Ch}}
\newcommand  {\Cokernel}{\operatorname{Coker}}
\newcommand  {\cosk}	{\operatorname{cosk}}
\newcommand  {\Cov} 	{{\rm\text{Cov}}}

\newcommand  {\Fil}	{\operatorname{Fil}}

\newcommand  {\gr}      {\operatorname{gr}}

\newcommand	{\hcH}[2]	{\overset{\hspace{-2pt}\raisebox{-1pt}[2pt][-2pt]{$\scriptscriptstyle{#1}$}}{\cH^{#2}}}
\newcommand  {\Hom}     {\operatorname{Hom}}

\newcommand  {\id}      {{\operatorname{id}}}
\newcommand  {\Image}   {\operatorname{Im}}

\newcommand  {\Kernel} 	{\operatorname{Ker}}

\newcommand  {\dirlim}  {\varinjlim}

\newcommand  {\lra}     {\longrightarrow}

\newcommand  {\op}      {\text{\rm op}}
\newcommand  {\ord}     {\operatorname{ord}}

\newcommand  {\pr}      {\operatorname{pr}}

\newcommand  {\quadand} {\quad\text{and}\quad}

\newcommand  {\ra}      {\rightarrow}

\newcommand  {\semi}    {\text{\rm semi}}
\newcommand  {\Semi}    {\operatorname{Semi}}

\newcommand  {\Set}     {{\text{\rm Set}}}

\newcommand  {\Sh}      {{\operatorname{Sh}}}

\newcommand  {\Tot}     {\operatorname{Tot}}

\newcommand{\uH}	{\underline{H}}

\newcommand{\RTC}	{ \text{\rm RTC}}

\def\mydate{\number\day\space\ifcase\month \or January\or February\or March\or 
April\or May\or June\or July\or
August\or September\or October\or November\or December\fi \space\number\year}

\begin{document}

\title[Rigidified torsor cocycles]
      {Rigidified torsor cocycles, hypercoverings and bundle gerbes}

\author[Stefan Schr\"oer]{Stefan Schr\"oer}
\address{Mathematisches Institut, Heinrich-Heine-Universit\"at,
40204 D\"usseldorf, Germany}
\curraddr{}
\email{schroeer@math.uni-duesseldorf.de}


\dedicatory{21 April 2019}

\begin{abstract}
We give a   geometric interpretation of sheaf cohomology for higher degrees $n\geq 1$ in terms of 
torsors on the member of degree  $d=n-1$ in hypercoverings of type $r=n-2$,  endowed with an additional data, the so-called rigidification.
This generalizes the   fact that   cohomology in degree one is the group of
isomorphism classes of torsors, where the rigidification becomes vacuous,
and that  cohomology in degree two can be expressed
in terms of bundle gerbes, where the rigidification becomes an associativity constraint.
\end{abstract}

\maketitle
\tableofcontents

\section*{Introduction}
\mylabel{Introduction}

Let $X$ be a topological space and $\shF$ be an abelian sheaf.
It is a classical fact that the first  sheaf cohomology $H^1(X,\shF)$, which is defined via injective resolutions,
has a   \emph{combinatorial description} in terms of  \v{C}ech cohomology $\cH^1(X,\shF)$, and also a \emph{geometric interpretation} as
the group $\pi_0(\shF\text{-Tors})$ of isomorphism
classes   for the  Picard category  of \emph{torsors}.
This was already discussed   by Grothendieck in \cite{Grothendieck 1955}, Section 5.1
and treated in utmost generality in Giraud's monograph  \cite{Giraud 1971}, Chapter III. 
In degree $n=2$, one may interpret cohomology via \emph{gerbes}, which are certain  fibered
categories.  

It is natural to ask whether    higher cohomology groups $H^n(X,\shF) $ also admit a   geometric interpretation.
The goal of this paper is to show that this indeed holds  in all degrees $n\geq 1$. 
The idea is to use torsors $\shT$ living on  pieces $U_{n-1}$  of  hypercoverings $U_\bullet$
of type $r=n-2$, endowed with some additional datum called \emph{rigidification}.
See Definition \ref{rtc} for details.
This generalizes Murray's notion of \emph{bundle gerbes}  \cite{Murray 1996}, which describe the situation  in
degree $n=2$ and are used in the context of differential geometry  and theoretical physics, mainly for the case
$\shF=\CC^\times$.
See   \cite{Murray 2010} for an introduction, \cite{Murray; Stevenson 2000} and \cite{Murray; Stevenson 2011} for further development,
and also  the work of Caray, Murray and Wang on 
\emph{higher bundle gerbes} \cite{Carey; Murray; Wang 1997}.

Our starting point is Lazard's observation that   higher cohomology groups  
already have a combinatorial description in terms of \emph{hyper-\v{C}ech cohomology} 
$$
H^n(X,\shF)=\hcH{r}{n}(X,\shF)=\dirlim H^n\Gamma(U_\bullet,\shF),
$$
where the direct limit runs over all hypercoverings $U_\bullet$ of type $r=n-1$.
Roughly speaking,   \emph{hypercoverings} are certain semi-simplicial coverings $U_\bullet$, which are more general that
\v{C}ech coverings, in the sense that one allows in finitely many degrees the passage to   open  coverings.
The \emph{type of a hypercoverings} indicates that
from a certain degree $r\geq 0$, one does not pass to open coverings anymore, such that $U_\bullet=\cosk_r(U_{\leq r})$
becomes a coskeleton. 
This beautiful idea was developed in utmost generality
in \cite{SGA 4b}, Expos\'e V, Section 7. 

In contrast to   common   \v{C}ech coverings,
where $U_n=U^{n+1}=U\times_X\ldots\times_XU$ are fiber products coming from a single over covering $U=\bigcup_{\lambda\in L} U_\lambda$,
the hypercoverings never gained widespread attention. Of course, this relies on the fact that
already for   paracompact spaces, sheaf cohomology coincides with \v{C}ech cohomology
(\cite{Godement 1964}, Chapter II, Theorem 5.10.1).
Such a result also holds for \'etale cohomology on quasicompact schemes admitting
an ample sheaf by Artin's result \cite{Artin 1971}, 
and for Nisnevich cohomology on quasicompact separated schemes by my findings
in  \cite{Schroeer 2017} and \cite{Schroeer 2019}.
In the realm of  homotopy theory, \v{C}ech coverings also suffice, because the canonical map $|U_\bullet|\ra X$ from the geometric realization
is a weak equivalence. This is due to Segal \cite{Segal 1968} for countable coverings,
and also holds in general as observed by Dugger and Isaksen \cite{Dugger; Isaksen 2004}.

However, hypercoverings   play a crucial role for the \'etale homotopy type, as introduced
by Artin and Mazur \cite{Artin; Mazur 1969}, and further studied in Friedlander's monograph \cite{Friedlander 1982}.
More generally, they appear in homotopical algebra, that is, abstract homotopy theory
in the context of  Quillen's model categories, as explained in Jardine's book \cite{Jardine 2015}.
Hypercoverings of type $r=1$ are essential to understand gerbes, as exposed by Giraud \cite{Giraud 1971},
and the subsequent differential geometry of gerbes, which was developed by Brylinski \cite{Brylinski 1993},
Hitchin \cite{Hitchin 1999} and Breen and Messing \cite{Breen; Messing 2005}, see also the survey \cite{Breen 2011}.
Despite all these advances, hypercoverings should become more popular, in my opinion.
Recently, I have used them to reduce Serre's Vanishing   for general affine schemes to the noetherian situation
\cite{Schroeer 2018}.

To obtain a geometric description for  higher cohomology groups $H^n(X,\shF)$, we introduce the notion of \emph{rigidified torsor cocycles}
$(U_\bullet,\shT,\varphi)$ for the abelian sheaf $\shF$ in degree $n\geq 1$. 
Here $U_\bullet$ is  a hypercovering of type $r=n-2$, and $\shT$ is a torsor for the abelian sheaf  $\shF|U_{n-1}$,
and the \emph{rigidification} $\varphi$ is a section over $U_n$ of the \emph{alternating preimage} 
$$
p_\alt^*(\shT)=p_0^*(\shT)\wedge p_1^*(\shT^{-1})\wedge\ldots\wedge p_n^*(\shT^{(-1)^n}),
$$
subject to the cocycle condition $q_\alt^*(\varphi)=0$. See Definition \ref{rtc} for details.
The $p_i:U_n\ra U_{n-1}$ and $q_j:U_{n+1}\ra U_n$ denote   face operators
in the hypercovering.
The rigidified torsor cocycles form a \emph{fibered Picard category} $\foR_{\shF}^n$, a notion going back
to MacLane \cite{MacLane 1963} and studied in depth by Deligne \cite{SGA 4c}, Expos\'e XVIII, Section 4.1.
The category $\foR_\shF^n$ contains  ``coboundary objects'' coming from alternating preimages of
$\shF|U_{n-2}$-torsors.  
Forming fiber-wise the resulting   residue class groups and passing to a filtered  direct limit along homotopy  classes of refinements, 
we obtain the abelian group $\operatorname{RTC}^n(\shF)$
of \emph{equivalence classes of rigidified torsor cocycles}.  Our main result is that this indeed gives a geometric description 
for sheaf cohomology:

\begin{maintheorem}(See Thm.~\ref{comparison bijective})
For each integer $n\geq 1$, there  is a natural identification of abelian groups $\RTC^n(\shF)=H^n(X,\shF)$.
\end{maintheorem}

To define the \emph{comparison map} $\RTC^n(\shF)\ra H^n(X,\shF)$, we analyze the spectral sequences
attached to the double complex $\Gamma(U_\bullet,\shI^\bullet)$, arising from the injective resolution $\shI^\bullet$ used to define sheaf cohomology
and the hypercovering $U_\bullet$.
This leads to the so-called \emph{three-term complex} $C^{n-1}\stackrel{\Psi}{\ra} C^n\stackrel{\Phi}{\ra} C^{n+1}$,
which captures a small but essential part of the double complex. The comparison map takes the form
$$
\RTC^n(\shF)\longmapsto \dirlim\Kernel(\Phi)/\Image(\Psi) = H^n(X,\shF).
$$
In degree $n=1$, everything boils down to the classical identification 
$$
H^1(X,\shF)=\cH^1(X,\shF)=\pi_0(\shF\text{-Tors});
$$
the hypercovering is constant, and the rigidification is vacuous.
In degree $n=2$, we recover Murray's notion of \emph{bundle gerbes}; the hypercovering reduces to a  \v{C}ech covering,
and the rigidification translates into an associativity constraint.
In some sense, our construction gives a higher-dimensional generalization of bundle gerbes.
The approach works for arbitrary sites $\foC$ satisfying     mild technical assumptions,
introduced mainly for the sake of exposition.
In this setting,  $X\in\foC$ denotes a final object.
One recurrent technical challenge is to understand the effect of simplicial homotopies between refinements of hypercoverings.

After the completion of the paper, I learned form Michel Murray, 
Daniel Stevenson and David Roberts that similiar results appear in the work of  
Duskin \cite{Duskin 1975}, Glenn \cite{Glenn 1982}, Beke \cite{Beke 2004} and Jardine \cite{Jardine 2009},
albeit in different forms and  in much more categorical settings.
One advantage of the present approach is that rigidified torsor cocycle yield a concrete way
to represent cohomology classes, and that we have an explicit comparison map to sheaf cohomology.
This could make the theory more accessible in fields outside category theory.

\medskip
The paper is organized as follows:
In Section \ref{Fibered   picard} we discuss the notion of fibered Picard categories $\foR\ra\foC$
and certain resulting direct limits.
Section \ref{Hypercoverings} contains basic facts on hypercoverings $U_\bullet$.
We use them to express sheaf cohomology with the cohomology
of the three-term complex $C^{n-1}\stackrel{\Psi}{\ra} C^n\stackrel{\Phi}{\ra} C^{n+1}$ 
in Section \ref{Three-term complex}, by using spectral sequences.
In Section \ref{Iterated alternating}, we introduce the notion of iterated alternating
preimages $p_\alt^*(\shF)$.
Section \ref{Rigidified torsor cocycles} contains the central definition of rigidified
torsor cocycles $(U_\bullet,\shF,\varphi)$. The main result appears in Section \ref{Comparison map}, where we 
define the comparison map $\RTC^n(\shF)\ra H^n(X,\shF)$ and establish its bijectivity.
In the closing Section \ref{Bundle gerbes} we   discuss   how Murray's bundle gerbes can be seen as
rigidified torsor cocycles.

\begin{acknowledgement}
I wish  to thank Michel Murray for informing me about the work of  Duskin, Beke and Jardine.
This research was conducted in the framework of the   research training group
\emph{GRK 2240: Algebro-geometric Methods in Algebra, Arithmetic and Topology}, which is funded
by the DFG. 
\end{acknowledgement}

\section{Fibered Picard categories}
\mylabel{Fibered   picard}

In this section, we   recall  some    notions from category theory that I found rather useful 
to phrase   results in later sections, 
and discuss certain direct limits attached to fibered Picard categories.

Let $\foR$ be a  \emph{symmetric monoidal category}. The monoidal structure
is given by a functor $\foR\times\foR\ra\foR$, which we usually write as $(A,B)\longmapsto A\wedge B$.
Moreover, we have  natural isomorphisms
$$
(A\wedge B)\wedge C\lra A\wedge (B\wedge C)\quadand A\wedge B\ra B\wedge A
$$
satisfying MacLane's axioms (see \cite{MacLane 1963} and \cite{SGA 4c}, Expos\'e XVIII, Section 4.1). 
These natural isomorphisms are   the \emph{associativity and commutativity constraints}.
If there is an   object $E\in\foR$ together with an natural isomorphism
$A\wedge E\ra A$, then the set $\pi_0(\foR)$ of isomorphism classes $[A]$, 
where $A\in \foR$ runs over all objects, acquires the structure of
an abelian monoid, with addition $[A]+[B]=[A\wedge B]$ and zero element $0=[E]$.

A  symmetric monoidal category $\foR$ is called \emph{strict} if 
the commutativity constraint for $A=B$ becomes the identity transformation $A\wedge A\ra A\wedge A$.
Note that this does not hold for the category $\foR=(R\text{-Mod})$ of modules over a ring, 
where the  monoidal structure is given by tensor products $M\otimes_RN$ and the commutativity constraint is $a\otimes b\mapsto b\otimes a$. 
But strictness actually holds on the full subcategory $\foR'=(R\text{-Inv})$ of invertible modules.

A \emph{groupoid} is a category in which all morphisms are isomorphisms.
A \emph{Picard category} is  a strictly symmetric monoidal groupoid $\foR$  
where the translation functors $X\mapsto A\wedge X$ are self-equivalences of categories, for all objects $A\in\foR$.
See \cite{Bertolin 2011} for a nice discussion.
The following is an important example: Let $X$ be a topological space,
$\shF$ be an abelian sheaf, and $\foR=(\shF\text{-Tors})$ be the Picard category of $\shF$-torsors $\shT$.
The monoidal structure is given by  $\shT'\wedge\shT''=\shF\backslash(\shT'\times\shT'')$, which is the quotient of
the product sheaf   by the diagonal  action.

For this category $\foR=(\shF\text{-Tors})$, we have a slight notational problem:
It is customary to write   group actions and torsor structures in a multiplicative way.
On the other hand, the group law for  abelian groups and abelian sheaves is preferable written additively, in particular
when it comes to \v{C}ech cohomology. So  the axiom  for a multiplicative action
$\shF\times\shT\ra\shT$, $(f,\varphi)\mapsto f\cdot \varphi$
of an additive group takes the form $(f+g)\cdot \varphi= f\cdot (g\cdot\varphi)$ and $0\cdot s=s$.
For $\shT=\shF$, we usually revert back to purely additive notation.

Next, we recall the notion of fibered categories.
For details, see \cite{SGA 1}, Expos\'e VI.
Given a functor  $F:\foR\ra\foC$ between arbitrary categories and an object $U\in\foC$, 
we write $\foR(U)$ for the \emph{fiber category}, consisting 
of objects and morphisms in $\foR$ send to the object $U$ and the morphism $\id_U$, respectively.
A morphism $f:A\ra B$ in $\foR$ over $\theta:U\ra V$ in $\foC$ is called \emph{cartesian}
if the map
$$
\Hom_U(A',A)\lra \Hom_\theta(A',B),\quad h\longmapsto f\circ h
$$
is bijective, for all $A'\in\foR(U)$.
One says the $\foR$ is an \emph{fibered category} if  for each morphism $\theta:U\ra V$
in $\foC$ and each object $B\in\foR(V)$ there is a  \emph{cartesian morphism} $f:A\ra B$ in $\foR$ inducing $\theta$,
and the composition of cartesian morphisms is cartesian.
If furthermore the fiber categories are groupoids, one says that $\foR$ is
a \emph{category fibered in groupoids}. Then all morphism
in $\foR$ are cartesian. The choice of a cartesian morphism
$f:A\ra B$, for each morphism $\theta:U\ra V$  and each object $B\in\shR(V)$,
is called a \emph{clivage}. One then regards $A$ as a ``fiber product'' $B\times_UV$ or ``base-change'' $\theta^*(B)$
of the object $B$ with respect to the morphism $\theta:U\ra V$.

Given furthermore a functor $\shR\times_\shC\shR\ra\shR$,
together with associativity and commutativity constraints,
that turn all fibers $\shR(U)$ into Picard categories,
we say that $\shR$ is a \emph{fibered Picard category}.
This notion was introduced by Deligne in \cite{SGA 4c}, Expos\'e XVIII, Section 1.4, where $\shC$ is furthermore
endowed with a Grothendieck topology and $\foR$ is regarded as  a \emph{Picard stack}.
In any case, we   get a contravariant functor 
\begin{equation}
\label{direct system fibers}
\foC\lra(\Ab),\quad U\longmapsto \pi_0(\foR(U))
\end{equation}
into the category of abelian groups.
The effect of a morphism $\theta:U\ra V$ is to send an isomorphism class
$[B]$ for the fiber category $\shR(V)$ to the isomorphism class   of its  base-change $A=\theta^*(B)=B\times_UV$.
Let us call $[B]\mapsto [A]$ the \emph{transition maps} for the
above functor.

In the category of sets, one now may form the \emph{direct limit} 
$$
\foR_\infty  = \dirlim \pi_0(\foR(U))
$$
for \eqref{direct system fibers}, regarded as a covariant functor on the opposite category $\foC^\op$.
Note that the group structures on the $\pi_0(\foR(U))$ do  not necessarily induce a group structure
on the direct limit, because direct limits do not necessarily commute with finite inverse limits.

For simplicity, suppose now that  $\foC$ admits a \emph{final object} $X\in\foC$, and that
the opposite category $\foC^\op$ is \emph{filtered}. 
For   the category  $\foC$, this  means that   each pair of morphisms $V'\ra V\leftarrow V''$  
sits in some commutative diagram
$$
\begin{CD}
U	@>>>	V''\\
@VVV		@VVV\\
V'	@>>>	V,
\end{CD}
$$
and for each pair of morphism $V\rightrightarrows V'$ there is a morphism
$U\ra V$ so that   the two compositions $U\rightrightarrows  V'$ coincide, confer \cite{SGA 4a}, Expos\'e I, Definition 2.7.
The first condition  is automatic if fiber products exist, and usually poses no problem in praxis, 
but the latter often becomes  tricky for categories
that do not come from ordered sets.
If both  conditions hold, the direct limit $\foR_\infty$ indeed inherits a group structure.
We call it  the \emph{stalk group} for the  fibered Picard category $\foR$
over the cofiltered category $\foC$. 

Often, the category $\foC$ is not cofiltered, but becomes cofiltered after passing to 
a \emph{quotient category} $\foC/{\sim}$, whose morphism are equivalence classes of morphisms in $\foC$,
and whose objects are the same.
In order to get a group structure on $\foR_\infty$, it then suffices to check that the 
the base-change map $\pi_0(\foR(V))\ra \pi_0(\foR(U))$   depends  only on the equivalence class of the morphism $\theta:U\ra V$.
In the application we have in mind, this independence actually takes place    after passing to
certain natural residue class groups $ \pi_0(\foR(U))/{\sim}$.

\section{Hypercoverings}
\mylabel{Hypercoverings}

In this section we recall some generalities on semi-simplicial coverings,
review the notion of hypercoverings and discuss the ensuing hyper-\v{C}ech cohomology. 
All results are basically contained in \cite{SGA 4b}, Expos\'e V, Section 2 and 7,
although in somewhat terse form. 
For generalities on  simplicial objects, we refer to the monographs
of May \cite{May 1967} and Weibel \cite{Weibel 1994}.

As customary, we write $\Delta$ for the category comprising the sets
$[n]=\{0,1,\ldots,n\}$ as objects and the monotonous maps $[m]\ra [n]$ as morphisms.
One may call this the \emph{category of simplex types}.
We will be only interested in the non-full subcategory $\Delta_\semi$,
which has the same object and where the morphisms are injective. 
The \emph{face map} $\partial_i:[n-1]\ra [n]$, $0\leq i\leq n$ is the unique injective monotonous map
that omits the value $i\in[n]$. These face maps satisfy the simplicial identities
$\partial_j\partial_i = \partial_i\partial_{j-1}$ for $i<j$.
A \emph{ semi-simplicial object} in some category $\foC$ is a contravariant functor
$U_\bullet:\Delta_\semi\ra\foC$. This amounts to  a sequence of objects $U_n\in\foC$, $n\geq 0$
together with \emph{face operators} $f_i=(\partial_i)^*:U_n\ra U_{n-1}$, $0\leq i\leq n$
satisfying the identities $f_if_j=f_{j-1}f_i$ for $i<j$.
Indeed, such data already specifies the entire semi-simplicial object in a unique way.
The semi-simplicial objects $U_\bullet$ form a category $\Semi(\foC)$, with natural transformations of functors as morphisms.
  
A \emph{homotopy} between two morphisms $\theta_\bullet,\zeta_\bullet:U_\bullet\ra V_\bullet$  
is a collection of morphisms $h_0,\ldots, h_n:U_n\ra V_{n+1}$,   given in each degree $n\geq 0$, satisfying 
$f_0h_0=\theta_n$ and $f_{n+1}h_n=\zeta_n$, together with  the relations
\begin{equation}
\label{homotopy}
f_ih_j=
\begin{cases}
h_{j-1}f_i 	& 	\text{if $0\leq i< j\leq n+1$;}\\
f_ih_{i-1}	&	\text{if $1\leq i=j\leq n$;}\\
f_ih_i	&	\text{if $1\leq i=j+1\leq n$;}\\
h_jf_{i-1}	&	\text{if $1\leq j+1<i\leq n+1$.}\\
\end{cases}
\end{equation}
We call them the \emph{simplicial homotopy identities}.
One   says that the two morphisms $\theta_\bullet$ and $\zeta_\bullet$ are \emph{homotopic} if they are equivalent under 
the equivalence relation generated by the homotopy relation. 
We then write $\Semi(\foC)/{\sim}$ for the ensuing \emph{quotient category},
where the objects are still the semi-simplicial objects in $\foC$, but the morphisms
are homotopy classes of natural transformations.
Note that the equivalence  relations generated by the homotopy relations 
on hom sets indeed from a \emph{congruence} in the sense of category theory (\cite{MacLane 1971}, Chapter II, Section 8).
 
Now suppose that $\foC$ is a \emph{site}, that is, a category endowed with a \emph{Grothendieck topology}
(confer \cite{SGA 4a}, Expos\'e II).
For simplicity, we assume that there is a \emph{final object} $X\in\foC$, and that
the topology comes from a \emph{pretopology}. The latter means that for each object $V\in\foC$, one has specified
a collection $\Cov(V)$ of \emph{covering families} $(U_\lambda\ra V)_{\lambda\in L}$ satisfying certain
axioms. For the sake of exposition, we also assume that for each covering family the disjoint
union $U=\bigcup_{\lambda\in L}U_\lambda$ exists. This ensures that one may refine each covering family to 
a one-member covering family, which one may call \emph{covering single}. 

One should have the following example in mind: Given a topological space $X$, let
$\foC$ be the category of  $X$-spaces  $(V,f)$ whose structure map $f:V\ra X$
is a local homeomorphism, and $\Cov(V)$ comprises  those families where  $\bigcup_{\lambda\in L} U_\lambda\ra V$ is surjective.
According to the Comparison Lemma (\cite{SGA 4a}, Expos\'e V, Theorem 4.1), the restriction functor $\Sh(\foC)\ra \Sh(X)$
is an equivalence of categories. In turn, sheaves and cohomology for the space  $X$ is the  same for the site $\foC$.

Let $\shG$ be an abelian presheaf on the site $\foC$.
For each semi-simplicial object $U_\bullet$ in $\foC$, we get the cochain complex
$\Gamma(U_\bullet, \shG)$.
The coboundary operators 
\begin{equation}
\label{horizontal}
\partial:\Gamma(U_n,\shG)\lra\Gamma(U_{n+1},\shG),\quad g\longmapsto\sum_{j=0}^{n+1} (-1)^jq_j^*(g) 
\end{equation}
are induced from the face operators $q_j:U_{n+1}\ra U_n$, $0\leq j\leq n+1$ (see for example \cite{Kashiwara; Schapira 2006}, Proposition 11.4.2).
This construction is  functorial with respect to the presheaf $\shG$ and the hypercovering $U_\bullet$.
Given two morphisms $\theta_\bullet,\zeta_\bullet:U_\bullet\ra V_\bullet$
and a   homotopy $h_0,\ldots, h_n:U_n\ra V_{n+1}$ between them, we can form the homomorphisms
$$
s:\Gamma(V_{n+1},\shG)\lra \Gamma(U_n,\shG), \quad g\longmapsto \sum_{i=0}^n (-1)^ih_i^*(g)
$$
The following fact appears, in one form or another, over and over in homological algebra.
The following form appears in  \cite{Stacks Project 2018}, Tag 019S, compare 
also the proof for Proposition \ref{homotopy torsors} below.

\begin{lemma}
\mylabel{homotopy cochains}
In the above situation, we have
$\theta^*_n-\zeta_n^* = s\partial - \partial s$ as homomorphisms $\Gamma(V_n,\shG)\ra\Gamma(U_n,\shG)$,
for all degrees $n\geq 0$.
\end{lemma}

In other words, simplicial homotopies yield cochain homotopies. In particular, the   homomorphisms $\theta_n$ and $\zeta_n$ 
induce the same map on cohomology.
More generally, the contravariant functor
$$
\Semi(\foC)\lra \Ch(\Ab),\quad U_\bullet\longmapsto \Gamma(U_\bullet,\shG)
$$
into the abelian category $\Ch(\Ab)$ of cochain complexes of abelian groups induces a functor
$$
\Semi(\foC)/{\sim}\,\lra K(\Ab)
$$
into the   quotient category $K(\Ab)=\Ch(\Ab)/{\sim}$, where the morphisms are  cochain maps modulo cochain homotopy.
Note that this is still an additive category, in fact a triangulated category, but no longer an abelian category.

A semi-simplicial object $U_\bullet $ in the site $\foC$ is called a \emph{semi-simplicial covering}
if  all face operators $p_i:U_n\ra U_{n-1}$ and the augmentation $U_0\ra X$ are coverings.
Morphisms between semi-simplicial coverings are called \emph{refinements}.
Using the covering families belonging to $\Cov(U_n)$, $n\geq 0$ we get categories 
of sheaves $\Sh(U_n)$. The face operators thus yield contravariant functors
$$
p_i^*:\Sh(U_{n-1})\lra \Sh(U_n), \quad\shG\longmapsto p_i^*(\shG).
$$
Here we write $p_i^*(\shG)=p_i^{-1}(\shG)$ for the preimage sheaf, because there is no danger of confusion.
In the presence of points and stalks, this means $p_i^*(\shG)_a = \shG_{p_i(a)}$. These preimage functors are related by  
natural isomorphisms $q_j^*\circ p_i^* \simeq p_i^*\circ p_{j-1}^*$ for $i<j$, which we regard
as an identifications.
Here 
$$
p_i:U_n\ra U_{n-1}\quadand q_j:U_{n+1}\ra U_n
$$
are the face operators defined in degree $n$
and degree $n+1$, respectively.  
In this context, we     usually write $r_k:U_{n-1}\ra U_{n-2}$ for   face operators in degree $n-1$.
This convention will be uses many times throughout.

Given a covering $U\ra X$, we may form the fiber products 
$$
U_n=U^{n+1}=U\times_X\ldots\times_XU,\quad n\geq 0
$$
and use   as face operators the projections $p_i:U_n\ra U_{n-1}$ that omit the $i$-th entry.
Such semi-simplicial coverings occur in the definition of \v{C}ech cohomology, and we refer to them
as \emph{\v{C}ech coverings}. Note that in degree one, the face operators become the projections $p_0=\pr_2$ and $p_1=\pr_1$.

\v{C}ech coverings are special cases of the   more flexible \emph{hypercoverings}, which  are certain semi-simplicial coverings $U_\bullet$   constructed 
by the following   recursive procedure:
One  formally starts with $U_{-1}=X$.
If for some degree $n\geq 0$  the objects $U_{-1},U_0,\ldots, U_{n-1}$ are already defined, 
one extends this \emph{truncated semi-simplicial} covering $U_{\leq n-1}$
to a full semi-simplicial covering $L_\bullet$   by taking fiber products in the universal way,
such that $\Hom(T_{\leq n-1},U_{\leq n-1})=\Hom(T_\bullet,L_\bullet)$ for all other
semi-simplicial coverings $T_\bullet$.
Now one is allowed to choose a covering family $(W_\lambda\ra L_n)_{\lambda\in L}$,
and defines $U_n=\bigcup_{\lambda\in L} W_\lambda$ as the corresponding covering single.
This concludes the recursive construction.
However, one demands that for some degree $r=n-1$, there is no passage to a covering family,
such that $U_\bullet=L_\bullet $. In turn, the hypercovering $U_\bullet$ is entirely determined
by the truncated covering $U_{\leq r}$. Equivalently, the morphism $U_\bullet\ra \cosk_r(U_{\leq r})$ 
into the \emph{coskeleton} is an isomorphism. One than says that $U_\bullet$ is a \emph{hypercovering of type $r\geq 0$}.
Note that the hypercoverings of type $r=0$ are precisely the   \v{C}ech coverings.
We refer to \cite{SGA 4b}, Expos\'e V, Section 7 for details.


Now let $\foH_{X,r}$ be the category of hypercoverings $U_\bullet $ of type $r\geq 0$.
Recall that the morphisms are called refinements. We   get a    functor 
$$
\foH_{X,r}^\op\lra\Ch(\Ab),\quad U\longmapsto \Gamma(U_\bullet,\shG)
$$
into the abelian category $\Ch(\Ab)$  of cochain complexes of abelian groups.
It satisfies Grothendieck's axiom (AB3), hence admits all direct sums and direct limits.
One may form the direct limit of the above functor, but looses control over
the resulting cochain groups, because the index category $\foH_{X,r}$ is in general not
filtered.

However, we may pass to the \emph{quotient category} $\overline{\foH}_{X,r}=\foH_{X,r}/{\sim}$,
where the morphisms are  homotopy classes of refinements.
Let us call  it the \emph{quotient category of hypercoverings}. 
Thus we formally get a functor
$$
\overline{\foH}^\op_{X,r}\lra K(\Ab),\quad  U\longmapsto \Gamma(U_\bullet,\shG).
$$
As explained in \cite{SGA 4b}, Expos\'e V, Theorem 7.3.2, the opposite quotient category $\overline{\foH}_{X,r}^\op$ becomes filtered.
In other words, the above may be regarded as an  \emph{ind-object} in the sense of \cite{SGA 1}, Expos\'e 1, Section 8.2.
Passing to cohomology, we   obtain    direct limits
$$
\hcH{r}{p}(X,\shG) = \dirlim  H^p\Gamma(U_\bullet,\shG)\in(\Ab).
$$
One may   form the   direct limits as  sets or     abelian groups, and use as index category
either $\foH_{X,r}^\op$ or the quotient category $\overline{\foH}^\op_{X,r}$, and always gets the same group.
These groups are called  \emph{hyper-\v{C}ech cohomology groups 
of type $r\geq 0$}, for each abelian presheaf $\shG$. For type $r=0$, this gives the  usual \v{C}ech cohomology $\cH^p(X,\shG)$.

Now write $(\text{AbP}/\foC)$ for the category of all abelian presheaves $\shG$. Generalizing from \v{C}ech cohomology,
one easily sees that for each fixed type $r\geq 0$, the 
$$
\hcH{r}{p}:(\text{AbP}/\foC)\lra(\text{Ab}),\quad\shG\longmapsto \hcH{r}{p}(X,\shG)
$$
form $\partial$-functors, and the canonical inclusions $\foH_{X,r}\subset\foH_{X,r'}$ for $r\leq r'$
induce natural transformations between $\partial$-functors.
For type $r=0$ we get the universal $\partial$-functor and the universal natural transformations.
Note further that the restrictions to the full subcategory
$(\text{Ab}/\foC)$ of all abelian sheaves $\shF$ usually fails to be a $\partial$-functor,
because the inclusion functor is not exact and the  subcategory contains more short exact sequences in general.
However, we have the following vanishing result:

\begin{proposition}
\mylabel{injective object acyclic}
Suppose  $\shI\in(\text{\rm AbP}/\foC)$   satisfies the sheaf axiom and   becomes an injective object in the abelian category $(\Ab/\foC)$.
Then for each hypercovering $U_\bullet$, the cochain complex $\Gamma(U_\bullet,\shI)$ has no cohomology in degrees $p\geq 1$.
In particular, $\shI$ is acyclic with respect to hyper-\v{C}ech cohomology for each type $r\geq 0$.
\end{proposition}

\proof
For each object $U\in\foC$, consider the corresponding representable presheaf $h_U:\foC\ra(\Set)$ given by $h_U(V)=\Hom_\foC(V,U)$.
The Yoneda Lemma yields   $\Hom_{(\text{\rm AbP}/\foC)}(h_U,\shG)=\Gamma(U,\shG)$ for each set-valued presheaf $\shG$.
Write $\ZZ h_U:\foC\ra(\Ab)$ for the resulting abelian presheaf, whose group of local section over $V$
is the free abelian group generated by the set $h_U(V)$.
We thus get $\Hom_{(\Ab/\foC)}(\ZZ h_U,\shG)=\Gamma(U,\shG)$, for each abelian presheaf $\shG$.
Finally, write $\ZZ_U$ for the sheafification of the abelian presheaf $\ZZ h_U$.
By the universal property of sheafification, we obtain $\Hom_{(\Ab/\foC)}(\ZZ_U,\shF)=\Gamma(U,\shF)$
for each abelian sheaf $\shF$.

Now let $U_\bullet $ be a hypercovering of type $r\geq 0$, and consider the resulting
semi-simplicial sheaf $\ZZ_{U_\bullet}$ and the ensuing chain complex of abelian sheaves
$$
\ldots\lra\ZZ_{U_2}\lra \ZZ_{U_1}\lra \ZZ_{U_0} 
$$
on $\foC$, where the boundary maps are alternating sums induced from the face operators $q_j:U_{n+1}\ra U_n$.
According to \cite{SGA 4b}, Expos\'e V, Theorem 7.3.2 this  chain complex is exact in degrees $p\geq 1$,
and the cokernel for the map on the right is $\ZZ_X$.
Since $\shI\in(\Ab/\foC)$ is an injective object, the functor $\shF\mapsto\Hom_{(\Ab/\foC)}(\shF,\shI)$ is exact.
Applying this exact hom functor to our  exact sequence and using the preceding paragraph, we get
an exact sequence
$$
0\lra \Gamma(X,\shI)\lra \Gamma(U_0,\shI)\lra \Gamma(U_1,\shI)\lra \ldots.
$$
Here the coboundary maps coincide with \eqref{horizontal}. In turn, the cohomology groups 
$H^p\Gamma(U_\bullet,\shI)$ vanish for $p\geq 1$.
Passing to direct limits with respect to $U_\bullet\in\foH_{X,r}^\op$, we get $\hcH{r}{p}(X,\shI)=0$
\qed

\section{The three-term complex}
\mylabel{Three-term complex}

We keep the notations and assumptions from the preceding section, such that $\foC$ is a site with
final object $X\in\foC$.
The goal now is to   express sheaf cohomology
in terms of certain three-term complexes, which arise from   total complexes and spectral sequences involving
hypercoverings. We use the standard conventions with regards to 
double complexes, compare \cite{Cartan; Eilenberg 1956}, Chapter IV and XV.

For each abelian sheaf $\shF$, choose once and for all an injective resolution
\begin{equation}
\label{vertical}
0\lra\shF\lra\shI^0\stackrel{d}{\lra}\shI^1\stackrel{d}{\lra}\ldots,
\end{equation}
such that $H^p(X,\shF)=H^p\Gamma(X,\shI^\bullet)$.
Furthermore, fix some hypercovering   $U_\bullet $. This gives    cochain complexes
$\Gamma(U_\bullet,\shI^q)$ with coboundary maps 
\begin{equation}
\label{horizontal again}
\partial:\Gamma(U_n,\shI^q)\lra\Gamma(U_{n+1},\shI^q),\quad f\longmapsto\sum_{j=0}^{n+1} (-1)^jq_j^*(f).
\end{equation}
The two types of coboundary maps are compatible, such that the  diagrams
$$
\begin{CD}
\Gamma(U_p,\shI^{q+1})	@>\partial_{p,q+1}>>	\Gamma(U_{p+1},\shI^{q+1})\\
@Ad_{p,q}AA						@AAd_{p+1,q}A\\
\Gamma(U_p,\shI^q)	@>>\partial_{p,q}> 	\Gamma(U_{p+1},\shI^q)
\end{CD}
$$
are commutative.
To obtain a \emph{double complex} $\Gamma(U_\bullet,\shI^\bullet)$, we use the usual sign trick and replace
the vertical differentials $d_{p,q}$  by $  (-1)^pd_{p,q}$.
Up to this sign change,  the vertical differential comes from the injective resolution \eqref{vertical}, whereas the horizontal differential
comes from the hypercovering \eqref{horizontal again}.

In turn, we obtain a \emph{total complex} $\Tot\Gamma(U_\bullet,\shI^\bullet)$.
First, we consider  the \emph{horizontal filtration}, which has 
$$
\Fil^q=\bigoplus_{p\geq 0, b\geq q} \Gamma(U_p,\shI^b).
$$
The associated graded complex is $\gr^q = \Gamma(U_\bullet,\shI^q)$, with differential
given by  \eqref{horizontal again}. This gives a spectral sequence with $E_1^{pq}= H^p\Gamma(U_\bullet,\shI^q)$.
Since the presheaf $\shI^q$ is an injective sheaf, we have 
$E_1^{pq}=0$ for all $p>0$, according to Proposition \ref{injective object acyclic}.
Furthermore, $E_1^{0,q}= H^0(X,\shI^q)$.
The differential $E_1^{0,q}\ra E_1^{0,q+1}$ is induced by \eqref{vertical},
and the definition of sheaf cohomology gives  $E_2^{0,q}=H^q(X,\shF)$.
Since there are no non-trivial differentials on the $E_2$-page, the spectral sequence
collapses, such that $E_2^{pq}=E_3^{pq}=\ldots=E_\infty^{pq}$.
In turn, the filtration on $\Fil^q H^{p+q}\Tot\Gamma(U_\bullet,\shF)$ has just one non-trivial step.
This shows:

\begin{proposition}
\mylabel{cohomology total complex}
For each hypercovering $U_\bullet $ and each degree $q\geq 0$, the edge map
$$
H^q\Tot\Gamma(U_\bullet,\shI^\bullet)\lra H^q(X,\shF)  
$$
for the above spectral  sequence is bijective.
\end{proposition}

The total complex is functorial in the hypercovering $U_\bullet$, and the edge map is compatible
with refinements. In turn,   refinements induces identities on the cohomology of the total space.

Next, consider   the \emph{vertical filtration} $\Fil^p=\bigoplus_{a\geq p, q\geq 0} \Gamma(U_a,\shI^q)$
on the total complex, whose associated graded complex is $\gr^p=\Gamma(U_p,\shI^\bullet)$.
This   gives  another  spectral sequence,
with   $E_1^{pq}= H^q\Gamma(U_p,\shI^\bullet)$ and $E_2^{pq}=H^p\Gamma(U_\bullet,\uH^q(\shF))$.
Here $\uH^q(\shF)$ denotes the presheaf defined by $\Gamma(V,\uH^q(\shF)) = H^q(V,\shF)$.
The spectral sequence is functorial in the hypercovering $U_\bullet$, and homotopic refinements induce identical maps.
Passing to the direct limit over all $U_\bullet\in\overline{\foH}_{X,r}^\op$ and using the preceding
proposition, we obtain  the  spectral sequence  
$$
E_2^{pq} = \hcH{r}{p}(X,\uH^q(\shF)) \Longrightarrow H^{p+q}(X,\shF), 
$$
as given in  \cite{SGA 4b}, Expos\'e V, Section  7.4.
According to loc.\ cit., Theorem 7.4.1 we have:

\begin{theorem}
\mylabel{edge map bijective}
For each degree $n\leq r+1$, the   edge map 
$$
\hcH{r}{n}(X,\shF)\lra H^n\Tot\Gamma(U_\bullet,\shI^\bullet) = H^n(X,\shF)
$$
for the above spectral sequence is bijective.
\end{theorem}

This should be seen as a far-reaching generalization of the fact that cohomology in degree one
coincides with \v{C}ech cohomology, which is the special case 
$r=0$ and $n=1$. 
Note, however, that in general the edge map is not bijective for $n=r+2$.

We shall exploit this failure as follows:  
Fix some integer $n\geq 1$, and consider  the three abelian groups
\begin{equation*}
\begin{split}
C^{n-1}	& = \Gamma(U_{n-2},\shI^0/\shF) \oplus\Gamma(U_{n-1},\shI^0),\\
C^{n}	&= \Gamma(U_{n-1},\shI^1)\oplus\Gamma(U_n,\shI^0), \\
C^{n+1}	&= \Gamma(U_{n-1},\shI^2)\oplus\Gamma(U_n,\shI^1)\oplus\Gamma(U_{n+1},\shI^0).	
\end{split}
\end{equation*}
Taking the horizontal and vertical differentials  from \eqref{horizontal again} and \eqref{vertical}, 
we get matrices 
$$
\Phi=\begin{pmatrix}
\epsilon d		&	0	\\
\partial	&	-\epsilon d	\\
0		&	\partial
\end{pmatrix}
\quadand
\Psi=\begin{pmatrix}
\partial	&	\epsilon d\\
0		& 	\partial
\end{pmatrix}
$$
that define a \emph{three-term cochain complex}
$C^{n-1}\stackrel{\Psi}{\ra} C^n\stackrel{\Phi}{\ra} C^{n+1}$. Here $\epsilon=(-1)^{n-1}$ is the sign introduced for
the double complex  $\Gamma(U_\bullet,\shI^\bullet)$.
One should have the following diagram in mind, which arises from  the double complex $\Gamma(U_\bullet,\shI^\bullet)$:
$$
\xymatrix{
						& \Gamma(U_{n-1},\shI^2)		&			\\
\Gamma(U_{n-2},\shI^0/\shF)\ar[r]^{\partial}	& \Gamma(U_{n-1},\shI^1)\ar[u]^{\epsilon d}\ar[r]^{\partial}	& \Gamma(U_n,\shI^1)\\
						& \Gamma(U_{n-1},\shI^0)\ar[u]^{\epsilon d}\ar[r]_{\partial}	& \Gamma(U_n,\shI^0)\ar[u]_{-\epsilon d}\ar[r]_{\partial}	& \Gamma(U_{n+1},\shI^0)
}
$$
Note that  upper left entry of the matrix $\Psi$ is given by the canonical inclusion $\shI^0/\shF\subset\shI^1$ coming from the differential
$ \shI^0\ra\shI^1$, which ensures that the upper left entry of $\Phi\circ\Psi$ vanishes.
The   exact sequence $0\ra \shI^0/\shF\ra\shI^1\ra\shI^2$ induces an exact sequence
$$
0\lra\Gamma(U_{n-2},\shI^0/\shF)\lra \Gamma(U_{n-2},\shI^1)\stackrel{-\epsilon d}{\lra} \Gamma(U_{n-2},\shI^2).
$$
In turn, we have a homomorphism from the three-term cochain complex to the total cochain complex, which induces
the map 
$$
\Kernel(\Phi)/\Image(\Psi)\lra H^n\Tot\Gamma(U_\bullet,\shI^\bullet).
$$
on cohomology.
Both sides depend  functorially on the hypercovering: Each refinement $U'_\bullet\ra U_\bullet$ 
gives a restriction map, and we thus get a contravariant functor 
$\foH_{X,r}\ra (\Ab)$ given by $U_\bullet\mapsto \Kernel(\Phi)/\Image(\Psi)$.

Suppose that $\theta_\bullet,\zeta_\bullet:U'_\bullet\ra U_\bullet$ are homotopic refinement.
Each of them induces a cochain map between the three-term complexes formed with  $U_\bullet$ and $U'_\bullet$.
Using that the operator $s$ in Lemma \ref{homotopy cochains} is natural in the sheaf $\shG$,
one easily infers that the induced maps $\Kernel(\Phi)/\Image(\Psi)\ra \Kernel(\Phi')/\Image(\Psi')$ on cohomology coincide. 

Passing to the direct limit and taking the identification in Proposition \ref{cohomology total complex} into account, we get 
a map 
$$
\dirlim \Kernel(\Phi)/\Image(\Psi)\lra \dirlim H^n\Tot\Gamma(U_\bullet,\shI^\bullet)=H^n(X,\shF),
$$
where the direct limit runs over all hypercoverings $U_\bullet\in\foH_{X,r}^\op$.

\begin{theorem}
\mylabel{cohomology three-term}
For each integer $n\geq 2$ and each type  $r\geq n-2$, the above map is bijective, such that we get an identification
$H^n(X,\shF)=\dirlim \Kernel(\Phi)/\Image(\Psi)$.
\end{theorem}

\proof
First, we check that the map in question is surjective.
Fix some cohomology class $[\alpha]\in H^n(X,\shF)$.
Represent the class on some suitable  hypercovering $U_\bullet$ of type $r$ by some tuple
$$
\alpha=(\alpha_0,\alpha_1,\ldots,\alpha_n),
$$
where each entry is a local section $\alpha_i\in \Gamma(U_i,\shI^{n-i})$ and the
tuple   is a cocycle in the total complex $\Tot\Gamma(U_\bullet,\shI^\bullet)$.
We now show by induction on $0\leq i\leq n-2$ that after passing to  finer hypercoverings,
the tuple becomes cohomologous to a tuple of the form
$\beta=(0,\ldots,0,\beta_{i},\ldots,\beta_n)$.
In the final case $i=n-2$, the remaining pair $(\beta_{n-1},\beta_n)$ gives a cocycle in the three-term complex
inducing the class $[\alpha]$.

This assertion is trivial for $i=0$. Now suppose that $1\leq i\leq n-3$, and that $\alpha_0=\ldots,\alpha_{i-1}=0$.  
Then $\alpha_i\in\Gamma(U_i,\shI^{n-i})$
vanishes in $\Gamma(U_i,\shI^{n-i+1})$, thus defines a class $[\alpha_i]\in H^{n-i}(U_i,\shF)$.
Choose some covering $W\ra U_i$ on which the cohomology class vanishes, in other words,
$\alpha_i|W$ lies in the image of $\Gamma(W,\shI^{n-i-1})$.
Now we use $i< n-2\leq r$, so we may refine the hypercovering $U_\bullet \in\foH_{X,r}$ of type $r$
and assume that $\alpha_i$ lies in the image of $\Gamma(U_i,\shI^{n-i-1})$.
Subtracting from $\alpha$ the ensuing coboundary, we achieve $\alpha_0=\ldots=\alpha_i=0$.
This shows that the canonical map is surjective.

It remains to check that our map is also injective. Suppose 
an $n$-cocycle of the form $\alpha=(0,\ldots,0,\alpha_{n-1},\alpha_n)$ is a coboundary
in the total complex $\Tot\Gamma(U_\bullet,\shI^\bullet)$. Then
there is an $(n-1)$-cochain $\beta=(\beta_0,\ldots,\beta_{n-1})$ mapping to $\alpha$.
Clearly, we may replace $\beta$ by some cohomologous cochain.
As in the preceding paragraph, we thus may assume that $\beta_0=\ldots=\beta_{n-3}=0$.
The entry $\beta_{n-2}\in\Gamma(U_{n-2},\shI^1)$ then vanishes in $\Gamma(U_{n-2},\shI^2)$,
thus can be regarded as an local section of $\shI^0/\shF$.
In turn, the pair $(\beta_{n-2},\beta_{n-1})$ is an element from $C^{n-1}$
mapping to $(\alpha_{n-1},\alpha_n)\in C^n$.
This means that the map in question is injective.
\qed

\section{Iterated alternating preimages}
\mylabel{Iterated alternating}

We keep the assumptions of the preceding section, such that $\foC$ is a site
with final object $X\in\foC$. Let $\shF$ be an abelian sheaf on $\foC$.
We now collect some observations how torsor behave in simplicial coverings.
Together with the  identification $H^n(X,\shF)=\dirlim \Kernel(\Phi)/\Image(\Psi)$ from
the previous section,
this will later lead to our geometric interpretation of sheaf cohomology.

First recall that one may view   $H^1(X,\shF)$ as the group of isomorphism
classes $\pi_0(\shF\text{-Tors})$ of $\shF$-torsors $\shT$.
We write the group law for $\shF$ additively and the group action on $\shT$ multiplicatively.
Let $\shT^{-1}$ be the same set-valued sheaf $\shT$, but endowed with
the inverse $\shF$-action, given on local sections by $(-f)\cdot s $.
If $\shF\ra\shF'$ is a homomorphism of abelian sheaves, we obtain an induced $\shF'$-torsor
$\shF'\wedge^\shF\shT$ as the quotient 
$$
\shF'\wedge^\shF\shT= \shF\backslash(\shF'\times\shT)
$$
by the diagonal $\shF$-action.
 Given  $\shF$-torsors $\shT_i$ indexed by a finite set $I$,
and a map $\epsilon:I\ra\{\pm 1\}$, the product sheaf $\prod_{i\in I}\shT_i^{\epsilon_i}$ becomes a torsor
under  the abelian product sheaf  $\shP=\prod_{i\in I}\shF_i$.
We define the \emph{contracted product} $\bigwedge_{i\in I}\shF_i$ as the induced $\shF$-torsor
$$
\begin{CD}
\bigwedge_{i\in I}\shT_i^{\epsilon_i} =  \shF\wedge^\shP (\prod_{i\in I}\shT_i^{\epsilon_i})
\end{CD}
$$
with respect to the addition map $\shP=\bigoplus_{i\in I}\shF\ra \shF$. In the additive group $H^1(X,\shF)=\pi_0(\shF\text{-Tors})$
of isomorphism classes, we then have
$$
[\bigwedge_{i\in I}\shT_i^{\epsilon_i}]= \sum_{i\in I} \epsilon_i [\shT_i].
$$
The sum vanishes if and only if the contracted product admits a global section.
Under suitable assumptions, the contracted product even acquires a \emph{canonical  section}.
The precise meaning of this locution becomes apparent in the proof for the following assertion:

\begin{lemma}
\mylabel{canonical section}
Suppose there is a free involution $\sigma:I\ra I$ on the index set 
such  that $\epsilon_{\sigma(i)} = -\epsilon_i$ and  $\shT_i=\shT_{\sigma(i)}$
for all $i\in I$. Then the contracted product $\bigwedge_{i\in I}\shT_i^{\epsilon_i}$ 
admits a canonical   section.
\end{lemma}

\proof
By taking contracted products of canonical sections, it suffices to treat the case that the index sets  $I$ contains merely two elements,
such that the contracted product takes the form $\shT\wedge\shT^{-1}$. 
This contracted product comes with a canonical map $\varphi:\shT\wedge\shT^{-1}\ra\shF$, which sends a pair of local sections
$(s,s')$ to the element $f$, defined via the condition $f\cdot s=s'$. Recall   that $\shT=\shT^{-1}$ as set-valued sheaves.
This map is well-defined, and equivariant with respect to the canonical $\shF$-actions.
It must be an isomorphism, because the category $(\shF\text{-tors})$ is a groupoid.
The preimage of the zero-section $0\in\Gamma(X,\shF)$
is the canonical global section of $\shT\wedge\shT^{-1}$.
\qed

\medskip
Now let $U_\bullet$ be any semi-simplicial covering.
Fix some integer $n\geq 0$, and consider the face operators $p_i:U_{n}\ra U_{n-1}$ for $0\leq i\leq n$.
Given an $\shF|U_{n-1}$-torsor $\shT$, we define the \emph{alternating preimage} as
$$
p_\alt^*(\shT) = p_0^*(\shT)\wedge p_1^*(\shT^{-1})\wedge\ldots\wedge p_n^*(\shT)^{(-1)^n},
$$
which is an $\shF|U_n$-torsor. This torsor may or may not be trivial.
However, we shall see that the  \emph{iterated alternating preimage}
$q_\alt^*(p_\alt^*(\shT))$ becomes trivial, and in fact
acquires a canonical section. Here  $q_j:U_{n+1}\ra U_n$ for $0\leq j\leq n+1$ denote face operators
defined in degree $n+1$. 
The reason is as follows: Set 
$$
\shT_{(j,i)}=q_j^*(p_i^*(\shT))\quadand \epsilon_{(j,i)}=(-1)^{i+j}.
$$
By definition, the iterated alternating preimage is the contracted product for the $\shF|U^{n+1}$-torsors
$\shT_{(j,i)}^{\epsilon_{(j,i)}}$, where  
the indices satisfy  $0\leq j\leq n+1$ and $0\leq i\leq n$.
The  simplicial identities $\partial_j\partial_i = \partial_i\partial_{j-1}$ for $i<j$
yield  canonical identifications $q_j^*(p_i^*(\shG))=q_i^*(p_{j-1}^*(\shG))$ that are natural in the sheaf $\shG$.
Consider the index set 
$$
L=[n]\times[n-1]=\{(j,i)\mid  \text{$0\leq i\leq n$ and $0\leq   j\leq n+1$}\},
$$
which has   cardinality $(n+1)(n+2)$, an even number. It comes with the partition
$$
L'=\{(j,i)\mid \text{$i<j$}\}\quadand L''=\{(j,i)\mid\text{$i\geq j$}\},
$$
and both of these sets have cardinality $(n+1)(n+2)/2$.
In fact, we have a canonical map
$$
\psi:L'\lra L'',\quad (j,i)\longmapsto (i,j-1),
$$
which reflects the simplicial identities. This map is injective, hence bijective, because
domain and range have the same cardinality. In turn, we get a canonical free involution $\sigma:L\ra L$,
with $\sigma|L'=\psi$ and $\sigma|L''=\psi^{-1}$. 
By construction, 
$$
\shT_{\sigma(j,i)}=\shT_{(j,i)}\quadand\epsilon_{\sigma(i,j)}=- \epsilon_{(i,j)}.
$$
So we may apply Lemma \ref{canonical section} to the iterated alternating preimage $q_\alt^*(p_\alt^*(\shT))$ 
and infer that it comes with a \emph{canonical section} $\varphi_\can\in\Gamma(U_{n+1},q_\alt^*(p_\alt^*(\shT))$,
which yields a  \emph{canonical identification}
$$
q_\alt^*(p_\alt^*(\shT))=\shF|U_{n+1}.
$$
This simple observation has   remarkable consequences: 

If $\shT$ is an $\shF|U_{n-1}$ torsor 
and $\varphi\in\Gamma(U_n,p_\alt^*(\shT))$ is a section of the alternating preimage of the torsor,
the   alternating preimage   of the section can be regarded as 
a section $q_\alt^*(\varphi)\in\Gamma(U_{n+1},\shF)$.
In turn, we get a ``cocycle condition'' $q_\alt^*(\varphi)=0$ for   pairs $(\shT,\varphi)$.

Furthermore, each $\shF|U_{n-2}$-torsor $\shT_{n-2}$ yields a ``coboundary pair'' $(\shT,\varphi)$,
where the $\shF|U_{n-1}$-torsor  $\shT=r_\alt^*(\shT_{n-2})$ is the alternating preimage, formed with
the face operators $r_k:U_{n-1}\ra U_{n-2}$,  
and $\varphi=\varphi_\can$ is the canonical section of the iterated
alternating preimage $p_\alt^*(\shT)=p_\alt^*(r_\alt^*(\shT_{n-2}))$.
Indeed,   coboundary pairs satisfy the cocycle condition:

\begin{proposition}
\mylabel{coboundary pair}
With the preceding notation, we have $q_\alt^*(\varphi_\can)=0$.
\end{proposition}

This proof for this innocuous assertion requires  a little preparation. Consider the triple index set
$$
M=\{(j,i,k)\mid \text{$0\leq j\leq n+1$ and $0\leq i\leq n$ and $0\leq k\leq n-1$} \}.
$$
This set has cardinality $(n+2)(n+1)n$, which is divisible by six.
In simplicial notation $M=[n+1]\times[n]\times[n-1]$.
Furthermore, the simplicial identities yield two
involutions $\sigma,\eta:M\ra M$, determined by the condition
$$
\sigma(j,i,k)=(i,j-1,k)\quadand \eta(j,i,k) = (j, k, i-1)
$$
for $i<j$ and $k<i$, respectively. Consider the dihedral permutation group $G=\langle\sigma,\eta\rangle$ inside the symmetric group $S_M$
generated by the involutions, which satisfy the relations $\sigma^2=\eta^2=e$.

\begin{lemma}
\mylabel{dihedral}
The dihedral  group $G$ has order $\ord(G)=6$, and the $G$-action on the index set $M$ is free.
\end{lemma}

\proof
Fix some element $(j,i,k)\in M$ with $k<i<j$. Its $G$-orbit consists of the following elements,
arranged in a hexagonal pattern:
$$
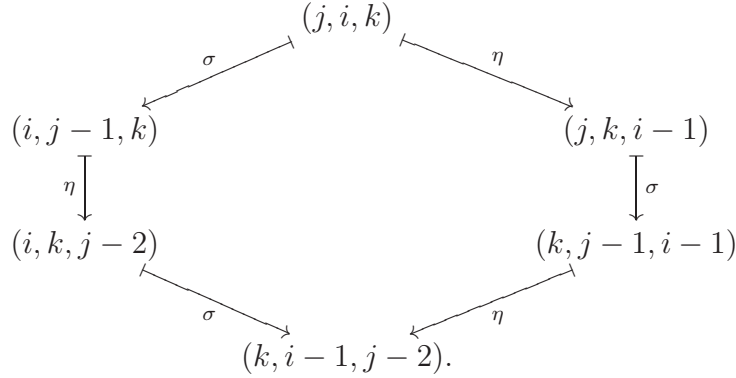
\xymatrix{
				&	(j,i,k)\ar@{|->}[dl]_{\sigma}\ar@{|->}[dr]^{\eta}\\
(i,j-1,k)\ar@{|->}[d]_{\eta}		&						& 	(j,k,i-1)\ar@{|->}[d]^{\sigma}\\
(i,k,j-2)\ar@{|->}[dr]_{\sigma}	&						& 	(k,j-1,i-1)\ar@{|->}[dl]^{\eta}\\
				&	(k,i-1,j-2).
}
$$
These six triples are pairwise different, because the entries $j,i,k$ are pairwise different.
Setting $g=\ord(G)$, we obtain  $6|g$.

By the simplicial identities, each element $(j',i',k')\in M$ lies in the  $G$-orbit of some element
$(j,i,k)$ with $k<i<j$. Hence all $G$-orbits comprise six elements, and we have the relation $\sigma\eta\sigma=\eta\sigma\eta$
in $G$. Therefore, the generators $\sigma,\eta\in G$ yield a surjection from
the dihedral group $D_3=\langle\sigma,\eta\mid \sigma^2=\eta^2=(\sigma\eta)^3\rangle$, thus $g|6$.
The upshot is that $g=6$, that the surjection $D_3\ra G$ is bijective, and that $G$ acts freely on $M$.
\qed

\medskip
\emph{Proof for Proposition \ref{coboundary pair}.}
By the very  definition, the iterated alternating preimage $p_\alt^*(\shT)=p_\alt^*(r_\alt^*(\shT_{U_{n-2}}))$
is a contracted product of the trivial $\shF|U_n$-torsors
\begin{equation}
\label{first torsor}
\shP_{ik} = p_i^*(r_k^*(\shT_{U_{n-2}}))^{(-1)^{i+k}}\wedge p_k^*(r_{i-1}^*(\shT_{U_{n-2}}))^{(-1)^{i+k-1}}
\end{equation}
for $k<i$, and the canonical section $\varphi_\can$ is the contracted product of   canonical sections $\varphi_{ik}\in\Gamma(U_n,\shP_{ik})$.
For each  $0\leq j\leq n+1$ with $i<j$, define
\begin{align}
\label{second torsor}
\notag
\shP_{j,i,k} 			&= q_j^*(\shP_{ik})^{(-1)^j}, \\
\shP_{k,j-1,i-1}        	&= q_k^*(\shP_{j-1,i-1})^{(-1)^k},\\
\notag
\shP_{i,j-1,k}			&= q_i^*(\shP_{j-1,k})^{(-1)^i}.
\end{align}
According to Lemma \ref{dihedral}, the   triple iterated alternating preimage  
$q_\alt^*(p_\alt^*(\shT))=q_\alt^*(p_\alt^*(r_\alt^*(\shT_{U_{n-2}})))$
is the contracted product of the trivial $\shF|U_{n+1}$-torsors
\begin{equation}
\label{third torsor}
\shP_{j,i,k} \wedge \shP_{k,j-1,i-1} \wedge \shP_{i,j-1,k},	
\end{equation}
where the contracted product runs over all triples $(j,i,k)$ with $k<i<j$. 
Substituting   \eqref{second torsor} and \eqref{first torsor}, we see that
\eqref{third torsor} is the contracted product of six torsors, which are   non-trivial in general.
To simplify notation,  set $\epsilon=(-1)^{i+j+k}$ and $\shS_{b,a,c}=q_b^*(p_a^*(r_c^*(((\shT_{U_{n-2}})))$, and 
arrange the six torsors in a hexagonal graph:
\begin{equation}
\label{hexagon}
\begin{gathered}
\xymatrix{
							& \shS_{j,i,k}^{\epsilon }\ar@{-}[dl]_{\sigma}\ar@{-}[dr]^{\eta}\\
\shS_{j,k,i-1}^{-\epsilon}\ar@{-}[d]_{\eta}		&									& \shS_{j,k,i-1}^{-\epsilon}\ar@{-}[d]^{\sigma}\\
\shS_{i,k,j-2}^{\epsilon }\ar@{-}[dr]_{\sigma}	&									& \shS_{k,j-1,i-1}^{\epsilon }\ar@{-}[dl]^{\eta}\\
							& \shS_{k,i-1,j-2}^{-\epsilon}
}
\end{gathered}
\end{equation}
Here the edges reflect the action of the generators $\sigma,\eta$ of  the dihedral permutation group $G\subset S_M$ on the set of indices $M$.
In light of the simplicial identities and the   exponents $\pm \epsilon$, any two adjacent torsors
are inverse to each other, in particular  have the same underlying set-valued sheaf.
The contracted product of the six torsors is a trivial $\shF|U_{n+1}$-torsor.
In fact, the contracted product has two \emph{canonical sections}, coming from the two possible cyclic  placement of parenthesis pairs:
Pairing the torsors incident with the $\sigma$-edges   gives 
$$
\left(\shS_{j,i,k}^{\epsilon }\wedge \shS_{j,k,i-1}^{-\epsilon} \right)\wedge 
\left(\shS_{i,k,j-2}^{\epsilon }\wedge \shS_{k,i-1,j-2}^{-\epsilon }\right)\wedge
\left(\shS_{k,j-1,i-1}^{\epsilon }\wedge\shS_{j,k,i-1}^{-\epsilon}\right),
$$
which yields the canonical section from the ``outer iterated preimage''. Pairing the torsors incident with the $\eta$-edges gives
$$
\left(\shS_{j,i,k}^{\epsilon   } \wedge \shS_{j,k,i-1}^{-\epsilon}\right)    \wedge
\left(\shS_{j,k,i-1}^{-\epsilon} \wedge \shS_{k,i-1,j-2}^{-\epsilon }\right) \wedge
\left(\shS_{i,k,j-2}^{\epsilon } \wedge \shS_{j,k,i-1}^{-\epsilon}\right),
$$
which yields the canonical section from the ``inner iterated preimage''.

Our task is to show that these two canonical sections coincide.
To verify this,   choose respective local sections 
$$
\xymatrix{
							& s_0\ar[dr]^{g_0}\ar@{|->}[dl]_{f_0}\\
s_1\ar@{|->}[d]_{g_1}		&									& s_5\ar@{|->}[d]^{f_1}\\
s_2\ar@{|->}[dr]_{f_2}	&									& s_4\ar@{|->}[dl]^{g_2}\\
							& s_3 
}
$$
for the torsors in the hexagonal diagram \eqref{hexagon}. The equations $s_1=f_0\cdot s_0$ and $s_2=g_1\cdot s_1$ etc.\ define
local sections $f_0,f_1,f_2$ and $g_0,g_1,g_2$ for the abelian sheaf  $\shF$, as indicated in the above diagram.
Obviously, they satisfy 
$f_0+g_1+f_2= g_0+f_1+g_2$,
hence $f_0-f_1+f_2=g_0-g_1+g_2$. The two sides of the latter equation correspond to  the canonical sections
of the sixfold contracted product, and our assertion follows.
\qed

\section{Rigidified torsor cocycles}
\mylabel{Rigidified torsor cocycles}

We keep the assumptions of the preceding section, such that $\foC$ is a site
with final object $X\in\foC$.
For simplicity, we assume  that the Grothendieck topology is given by a pretopology
$\Cov(V)$, $V\in\foC$ of covering families $(U_\lambda\ra V)_{\lambda\in L}$.
Furthermore, we suppose that for each covering family $(U_\lambda\ra V)_{\lambda\in L}$,
the disjoint union $U=\bigcup U_\lambda$ exists. 
We now introduce the central notion of this paper:

\begin{definition}
\mylabel{rtc}
Let $\shF$ be an abelian sheaf on the site $\foC$, and $n\geq 1$ be some integer.
A \emph{rigidified torsor cocycle} for $\shF$ in degree $n$  is a triple $(U_\bullet,\shT,\varphi)$
consisting of the following data:
\begin{enumerate}
\item $U_\bullet$ is a hypercovering of type $r=n-2$;
\item $\shT$ is a torsor for the abelian sheaf $\shF|U_{n-1}$  over $U_{n-1}$;
\item $\varphi$ is a   section of the alternating preimage $p_\alt^*(\shT)$  over $U_n$.
\end{enumerate}
These data satisfy   $q^*_\alt(\varphi)=0$
as section over $U_{n+1}$ of the iterated alternating preimage  
under the canonical identification $q_\alt^*(p_\alt^*(\shT))=\shF|U_{n+1}$.
\end{definition}

For simplicity, we also say that $(U_\bullet, \shT,\varphi)$ is a \emph{rigidified $\shF$-torsor $n$-cocycle}.
Recall that $p_i:U_n\ra U_{n-1}$ and $q_j:U_{n+1}\ra U_n$ denote  the   face operators
in the semi-simplicial covering $U_\bullet$ defined in degree $n$ and $n+1$, respectively.
One should have in mind the special case that the hypercovering $U_\bullet$ has type $r=0$,
such that the $U_m=U^{m+1}=U\times_X\times\ldots\times_XU$ form a  \v{C}ech covering.
Also note that a hypercovering of type $r=-1$ is just the constant semi-simplicial covering
with $U_m=X$ for all $m\geq 0$.

Some rigidified torsor cocycles arise as follows: Given a hypercovering $U_\bullet$ of type $n-2$,
an  $\shF|U_{n-2}$-torsor $\shT_{U_{n-2}}$ and a local section $s\in\Gamma(U_{n-1},\shF)$, 
we may form the alternating preimage $\shT=r_\alt^*(\shT_{U_{n-2}})$,
where the  $r_k:U_{n-1}\ra U_n$, $0\leq k\leq n-1$ denote   face operators defined in degree $n-1$.
Being an iterated alternating preimage,
the $\shF|U_n$-torsor $p_\alt^*(\shT)=p_\alt^*(r_\alt^*(\shT_{U_{n-2}}))$ on $U_n$ comes with a canonical section $\varphi=\varphi_\can$,
which we can multiply with $p^*_\alt(s)\in\Gamma(U_{n-1},\shF)$.
It follows from Proposition \ref{coboundary pair} that the triple $(U_\bullet,r_\alt^*(\shT_{U_{n-2}}),p^*_\alt(s)\cdot\varphi_\can)$ 
is a rigidified torsor cocycle. Let us call them \emph{rigidified torsor coboundaries}.
 
The rigidified $\shF$-torsor $n$-cocycles form a category $\foR_\shF^n$, 
where the morphisms 
$$
(\theta_\bullet,\tau): (U_\bullet, \shT,\varphi) \lra (U'_\bullet, \shT',\varphi')
$$
consists of a refinement $\theta_\bullet:U_\bullet\ra U'_\bullet$, together with
a morphism  $\tau: \shT\ra \shT'|U_{n-1}$ of $\shF|U_{n-1}$-torsors, such that $p_\alt^*(\tau)(\varphi)=\varphi'|U_n$.
Recall that $\foH_{X,n-2}$ denotes the category of hypercoverings $U_\bullet$ of type $r=n-2$,
and consider the forgetful functor
$$
\foR_\shF^n\lra\foH_{X,n-2},\quad (U_\bullet, \shT,\varphi)\lra U_\bullet.
$$
Clearly, all fiber categories $\foR_\shF^n(U_\bullet)$ are groupoids,
and all morphisms in $\foR_\shF^n$ are cartesian.
The fiber categories $\foR_\shF^n(U_\bullet)$ are actually    Picard categories.
The monoidal structure is given by 
$$
(U_\bullet,\shT,\varphi) = (U_\bullet,\shT',\varphi')\wedge(U_\bullet,\shT'',\varphi''),
$$
where  $\shT=\shT'\wedge\shT''$ and $\varphi=\varphi'\wedge\varphi''$. The associativity and commutativity
constraints come from the corresponding constraints for torsors.
The neutral object is given by the triple $(U_\bullet, \shF|U_{n-1}, 0)$.
The monoidal structure on   fibers comes from  the functor
$\foR_\shF^n\times_{\foH_{X,n-2}}\foR_\shF^n\ra \foR_\shF^n$ between categories over $\foH_{X,n-2}$,
where a triple consisting of $(U'_\bullet,\shT',\varphi')$ and $(U''_\bullet,\shT'',\varphi'')$ together
with an isomorphism $\theta_\bullet:U'_\bullet\ra U''_\bullet$ is mapped to the obvious wedge product taking $\theta_\bullet$ into account.
Summing up:

\begin{proposition}
\mylabel{forgetful functor}
The forgetful functor $\foR_\shF^n\ra\foH_{X,n-2}$ is a   fibered Picard category.
\end{proposition}

To proceed, let $\foT_\shF^{n-1}$ be the fibered Picard category whose objects are triples $(U_\bullet,\shT_{U_{n-2}},s)$, where $U_\bullet$
is a hypercovering of type $n-2$, and  $\shT_{U_{n-2}}$ is an $\shF|U_{n-2}$-torsor, and $s\in\Gamma(U_{n-1},\shF)$ is  a local section.
Morphism 
$$
(U_\bullet,\shT_{U_{n-2}},s)\lra (U'_\bullet,\shT'_{U_{n-2}},s')
$$
consists of a refinement
$\theta_\bullet:U_\bullet\ra U'_\bullet$ such that $s=s'|U_{n-1}$, 
together with a morphism  $\tau:\shT_{U_{n-2}}\ra \shT'_{U_{n-2}}|U_{n-2}$ of $\shF|U_{n-2}$-torsor.

Note that in the special case $n=1$, our hypercovering $U_\bullet$ of type $r=-1$ takes constant values $U_d=X$, $d\geq 0$.
Hence $\shT_{U_{-1}}$ is a torsor over $X$ endowed with a global section $s\in\Gamma(X,\shT)$.
In particular, all objects in the fiber categories $\foT^{-1}_{\shF}(U_\bullet)$ are isomorphic.

We have a commutative diagram  
$$
\xymatrix{
\foT_\shF^{n-1}\ar[dr]\ar[rr]	&		& \foR_\shF^n\ar[dl]\\
				& \foH_{X,n-2},	
}
$$
where the horizontal functor sends the triple  $(U_\bullet,\shT_{U_{n-2}},s)$ to the rigidified torsor coboundary
$(U_\bullet, r^*_\alt(\shT_{U_{n-2}}),p^*_\alt(s)\cdot\varphi_\can)$, and the two diagonal arrows  are the forgetful functors.
The horizontal functor respects the monoidal structure.
The cokernels
$$
\Cokernel\left(\pi_0(\foT_\shF^{n-1}(U_\bullet))\lra \pi_0(\foR_\shF^n(U_\bullet))\right)
$$
for the resulting homomorphisms between fiber-wise groups of isomorphism classes depend  functorially on the
hypercovering $U_\bullet$. The transition maps arise from base-changes.
We write the ensuing contravariant functor as
$$
\foH_{X,n-2}\lra (\Ab),\quad U_\bullet\longmapsto \pi_0(\foR_\shF^n(U_\bullet)) / \pi_0(\foT_\shF^{n-1}(U_\bullet)).
$$
The ensuing direct limit 
$$
\operatorname{RTC}^n(\shF) = \dirlim \pi_0(\foR_\shF^n(U_\bullet)) / \pi_0(\foT_\shF^{n-1}(U_\bullet))
$$
in the category of sets is called the set of \emph{equivalence classes of rigidified torsor cocycles} for the abelian sheaf $\shF$ in degree $n\geq 1$.
Here the  direct limit runs over all objects $U_\bullet\in\foH_{X,n-2}^\op$.

It is not a priori clear that the set $\RTC^n(\shF)$ inherits a group structure, because the category $\foH_{X,n-2}^\op$
is usually not filtered. 
We thus have to check that the transition maps depend only on homotopy classes of refinements,
that is, the functor factors over the quotient category  $\overline{\foH}_{X,n-2}^\op$.
To see this, suppose we have a homotopy $h_0,\ldots,h_n:U_n\ra V_{n+1}$ between two refinements $\theta_\bullet,\zeta_\bullet:U_\bullet\ra V_\bullet$,
as discussed in Section \ref{Hypercoverings}.    Write
$h_\alt^*(\shP)= h_0^*(\shP)\wedge h_1^*(\shP^{-1})\wedge \ldots\wedge h_n^*(\shP^{(-1)^n})$ for each $\shF|V_{n+1}$-torsor $\shP$.
 
\begin{lemma}
\mylabel{homotopy torsors}
Let $V_\bullet$ be a semi-simplicial covering and $\shT$ be an $\shF|V_{n-1}$-torsor.
With the above notation, we have a canonical  identification   
\begin{equation}
\label{identification with homotopy}
\theta_{n-1}^*(\shT)\wedge  r_\alt^*(h_\alt^*(\shT^{-1}))= \zeta^*_{n-1}(\shT) \wedge   h^*_\alt(p_\alt^*(\shT) )
\end{equation}
of $\shF|U_{n-1}$-torsors, which is natural in $\shT$. Here $p_i$ and $r_k$ denote face operators defined in degree $n$ and $n-1$, respectively.
\end{lemma}

\proof
The torsor $h^*_\alt(p_\alt^*(\shT))$ on the right hand side is a wedge product of torsors
$h_j^*(p_i^*(\shT))^\epsilon$, with exponents $\epsilon=(-1)^{i+j}$ and      indices $0\leq i\leq n$ and $0\leq j\leq n-1$.
Now recall from Section \ref{Hypercoverings} the  simplical homotopy identities   \eqref{homotopy}.
In the two boundary cases $i=j=0$ and $i=j+1=n$ we get $\theta_{n-1}^*(\shT)$ and $\zeta^*_{n-1}(\shT)^{-1}$, respectively.
The two  diagonal strips $1\leq i=j\leq n$ and $1\leq i=j+1\leq n+1$ give torsors  $h_i^*(p_i^*(\shT))$
and $h_i^*(p_{i-1}^*(\shT))^{-1} $ that are inverse to each other.
The remaining $n(n+1)-2-2(n-1)=n(n-1)$ torsors take the form
$$
h_j^*(p_i^*(\shT))^{\epsilon}=
\begin{cases}
r_i^*(h_{j-1}^*(\shT))^{\epsilon} 	& 	\text{if $0\leq i< j\leq n+1$;}\\
r_{i-1}^*(h_j^*(\shT))^{\epsilon}	&	\text{if $n+1\geq i>j+1>0$.}\\
\end{cases}
$$
These comprise the factors in the torsor $ r_\alt^*(h_\alt^*(\shT^{-1}))$ occurring on the left hand side of the natural identification.
Summing up, the simplicial homotopy identities, together with $\shP\wedge\shP^{-1}=\shF|U_{n-1}$ give the  desired  identification, which  is 
therefore natural in $\shT$.
\qed

\medskip
This ensures the desired homotopy invariance:

\begin{proposition}
\mylabel{transition maps}
Let $(V_\bullet, \shT,\varphi)$ be a rigidified torsor cocycle in degree $n\geq 1$,
and   $\theta_\bullet,\zeta_\bullet:U_\bullet\ra V_\bullet$ be two homotopic refinements.
Then the resulting   transition   maps
$$
\pi_0(\foR_\shF^n(V_\bullet))/ \pi_0(\foT_\shF^{n-1}(V_\bullet)) \lra \pi_0(\foR_\shF^n(U_\bullet))/ \pi_0(\foT_\shF^{n-1}(U_\bullet))
$$
induced by $\theta_\bullet$ and $\zeta_\bullet$ coincide.
\end{proposition}

\proof
It suffices to treat the case that there is a homotopy $h_0,\ldots,h_d:U_d\ra V_{d+1}$ between the refinements 
$\theta_\bullet$ and $\zeta_\bullet$. The idea is to extend the natural identification of $\shF|U_{n-1}$-torsors
in Lemma \ref{homotopy torsors} to an isomorphism of rigidified torsor cocycles, where the additional factors
$$
\shT'= r_\alt^*(h_\alt^*(\shT))\quadand\shT''= h^*_\alt(p_\alt^*(\shT) )
$$
are equipped with the structure of rigidified torsor coboundaries.
 
For $\shT'$, we take the canonical section $\varphi'=\varphi_\can$ for the iterated alternating preimage 
$q_\alt^*(\shT')$ of $h_\alt^*(\shT)$.
Note that in the special case $\shT=\shF|V_{n-1}$, this becomes $\shT'=\shF|U_{n-1}$ and $\varphi'=0$.

For $\shT''$, we use the given   $\varphi\in\Gamma(V_n,p^*_\alt(\shT))$
to define $\varphi''=-p_\alt^*(h_\alt^*(\varphi))$.
The ensuing triple $(U_\bullet,\shT'',\varphi'')$ is indeed a rigidified coboundary:
Using the identification $\shT''=\shF|U_{n-1}$ stemming from $h_\alt^*(\varphi)$,
we get $\shT''=\shF|U_{n-1}=r^*_\alt(\shF|U_{n-2})$ and $\varphi''=0$, and in particular $q^*(\varphi'')=0$.

It remains to see that the natural identification \eqref{identification with homotopy} is compatible with
the   four rigidifications $\theta_n^*(\varphi), \varphi',\zeta_n^*(\varphi),\varphi''$.
Applying $p^*_\alt$   and using that $\theta_\bullet$ and $\zeta_\bullet$ are natural transformations,
we get
$$
\theta_n^*(p_\alt^*(\shT))\wedge p^*_\alt(\shT') =
\zeta_n^*(p_\alt^*(\shT))\wedge p_\alt^*(\shT'').
$$
This is an identification of $\shF|U_n$-torsors, which is natural in the $\shF|V_{n-1}$-torsor $\shT$.
Our task is to verify   
\begin{equation}
\label{rigidifications coincide}
\theta_n^*(\varphi)\wedge \varphi'=\zeta_n^*(\varphi)\wedge\varphi'',
\end{equation}
as sections in the above torsor. 
This problem is local in $U_n$; to check it we may refine our hypercoverings, and even   increase their type. 
Now choose some covering $\tilde{V}_{n-1}\ra V_{n-1}$ on which $\shT$ acquires a section,
and form the fiber product $\tilde{U}_{n-1}=\tilde{V}_{n-1}\times_{V_{n-1}} U_{n-1}$.
Setting $\tilde{V}_{d}=V_d$ and $\tilde{U}_d=U_d$ for $d<n-1$, we obtain
truncated semi-simplicial coverings and   hypercoverings 
$$
\tilde{U}_\bullet=\cosk_{n-1}(\tilde{U}_{\leq n-1})\quadand
\tilde{V}_\bullet=\cosk_{n-1}(\tilde{V}_{\leq n-1})
$$
of type $r=n-1$. In order to check \eqref{rigidifications coincide}, we may restrict 
to the new hypercoverings. Replacing the old hypercoverings by the new ones and using the naturality of our construction,
we   have reduced the problem to the special case $\shT=\shF|V_{n-1}$ and thus  $\varphi'=0=h_\alt^*(q_\alt^*(\varphi))$. 

Consider the cochain complexes  $\Gamma(V_\bullet,\shF)$ and $\Gamma(U_\bullet,\shF)$ of abelian groups,
where the coboundary operators $\partial$ are   alternating sums as in \eqref{horizontal}.
The refinements induce cochain maps  $\theta_\bullet^*,\zeta_\bullet^*:\Gamma(V_\bullet,\shF)\ra\Gamma(U_\bullet,\shF)$.
With the notation from Lemma \ref{homotopy cochains}, we have
$
\theta_d^*-\zeta^*_d=  s\partial-\partial s
$
in all degrees $d\geq 0$. For $d=n$, we find
\begin{equation}
\label{reformulation coincide}
\theta_n^*(\varphi) +\varphi' = \theta_n^*(\varphi) - h_\alt^*(q_\alt^*(\varphi))= \zeta_n^*(\varphi)-p_\alt^*(h_\alt^*(\varphi))=  \zeta_n^*(\varphi)+\varphi''.
\end{equation}
It follows that  the equality \eqref{rigidifications coincide} holds.
\qed

\medskip
Let us unravel our set-up in degree $n=1$: Recall that a hypercovering $U_\bullet$ of type $r=-1$
is given by $U_m=X$ for all $m\geq 0$, and every face operator is the identity on $X$.
Hence the alternating preimage for the face operators $p_0,p_1:U_1\ra U_0$
is given by $p_\alt^*(\shT)=\shT\wedge\shT^{-1}=\shF$.
Likewise, the iterated alternating preimage is 
$$
q_\alt^*(p_\alt^*(\shT))=q_\alt^*(\shF)=\shF\wedge\shF^{-1}\wedge\shF=\shF,
$$
such that
$q_\alt^*(\varphi)=\varphi\varphi^{-1}\varphi=\varphi$.
In turn, a rigidified torsor cocycle $(U_\bullet,\shT,\varphi)$  in degree $n=1$
is entirely determined by the $\shF$-torsor $\shT$ over the final object $U_{0}=X$,
and $\varphi$ is the zero-section of $\shT\wedge\shT^{-1}=\shF|U_n$.
Using the description of $H^1(X,\shF)=\pi_0(\shF\text{-Tors})$, we immediately get:

\begin{proposition}
\mylabel{torsors}
There is a functorial identification $\RTC^1(\shF)=H^1(X,\shF)$.
\end{proposition}

The situation is more challenging in degree $n=2$:
Suppose we have an $\shF$-gerbe $\foG\ra\foC$.
Choose a covering $U\ra X$ over which there is an object $T\in\foG_U$.
The fiber produces $U_m=U^{m+1}$ form a \v{C}ech covering $U_\bullet $, that is, a hypercovering  of type $r=0$.
Write
$$
\xymatrix{
U_3\ar@<3pt>[r]^{q_j}	\ar@<1pt>[r]	\ar@<-1pt>[r]\ar@<-3pt>[r]	&	U_2\ar@<2pt>[r]^{p_i}\ar[r]\ar@<-2pt>[r]	& 	U_1\ar@<1pt>[r]^{r_k}\ar@<-1pt>[r]	&	U_0\ar[r]	&	X
}
$$
for the face operators defined in degrees $\leq 3$.
We get two objects $T_0=r_0^*(T)$ and $T_1=r_1^*(T)$ in the fiber category $\foG_{U_1}$, and thus an  $\shF|U_1$-torsor $\shT=\underline{\Hom}(T_1,T_0)$. 
To compute the  alternating preimage $p_\alt^*(\shT)$, we use the three simplicial identities 
$r_0p_1=p_0r_1$ and  $r_0p2=r_1p0$ and $r_1p_2=r_1p_1$ to get 
\begin{gather*}
p_0^*(\shT)^{-1} =  \underline{\Hom}(p_0^*r_0^*T,p_0^*r_1^*T),\\
p_1^*(\shT)^{-1}  =  \underline{\Hom}(p_1^*r_0^*T, p_1^*r_1^*T) = \underline{\Hom}(p_0^*r_0^*T, p_1^*r_1^*T),\\
p_2^*(\shT)^{-1}  =  \underline{\Hom}(p_2^*r_0^*T, p_2^*r_1^*T) = \underline{\Hom}(p_0^*r_1^*T, p_1^*r_1^*T).
\end{gather*}
Composition for  hom sets yields a canonical isomorphism of $\shF|U^2$-torsors
\begin{equation}
\label{hom sets}
\varphi: p_2^*(\shT)^{-1}\wedge p_0^*(\shT)^{-1}   \lra p_1^*(\shT)^{-1}, \quad
(h_2,h_0)\longmapsto h_2\circ h_0,
\end{equation}
which can be regarded as a section $\varphi\in\Gamma(U_2,p_\alt^*(\shT))$.
This gives us a triple $(U_\bullet,\shT,\varphi)$. 

\begin{proposition}
\mylabel{gerbes}
The triple  $(U_\bullet,\shT,\varphi)$ is a rigidified torsor cocycle.
\end{proposition}

\proof
We have to check that $q_\alt^*(\varphi)=0$ as section of the iterated alternating preimages $q_\alt^*(p_\alt^*(\shT))=\shF|U^3$.
Choose a covering $V\ra U_1$ over which the $\shF|U_2$-torsor $\shT$ acquires a section.
This section defines an isomorphism $r_1^*(T)|V\ra r_0^*(T)|V$ of objects in
the fiber category $\foG_V$, and we write 
 $h:r_0^*(T)|V\ra r_1^*(T)|V$ for its inverse.
Consider the truncated semi-simplicial covering  $V_{\leq 2}$
with $V_2=V$ and $V_1=U_1$ and $V_0=U_0$, and the ensuing hypercovering $V_\bullet=\cosk_2(V_{\leq 2})$ of type two.
The composition   \eqref{hom sets} of hom sets
defines via the equation
$$
p_2^*(h) \circ p_0^*(h) = f\cdot p_1^*(h)
$$
a   local section $f\in\Gamma(V_2,\shF)$.
Using the six simplicial identities $p_iq_j=p_{j-1}q_i$ for $i<j$,
the above equality gives four equations
\begin{gather*}
q_0^*p_2^*(h)\circ q_0^*p_0^*(h) = q_0^*(f)\cdot q_0^*p_1^*(h),\\
q_1^*p_2^*(h)\circ q_0^*p_0^*(h) = q_1^*(f)\cdot q_1^*p_1^*(h),\\
q_2^*p_2^*(h)\circ q_0^*p_1^*(h) = q_2^*(f)\cdot q_1^*p_1^*(h),\\
q_2^*p_2^*(h)\circ q_0^*p_2^*(h) = q_3^*(f)\cdot q_1^*p_2^*(h).
\end{gather*}
Here $q_j^*(f)\in\Gamma(V_3,\shF)$, and the  equations hold as morphisms $q_j^*p_i^*r_0^*(T)\ra q_j^*p_i^*r_1^*(T)$ in the fiber category
$\foG_{V_3}$.
In the above four equations, each of the terms $q_j^*p_i^*(h)$ 
appears twice, and the pair members come with ``opposite signs'' in the alternating preimage. This ensures $q_0^*(f)+q_2^*(f)=q_1^*(f)+q_2^*(f)$.
By definition of the iterated alternating preimage, one has 
$$
q_\alt^*(\varphi)|V_3= q_0^*(f)-q_1^*(f)+ q_2^*(f)-q_3^*(f).
$$
In turn, we have $q_\alt^*(\varphi)|V_3=0$. Since the induced morphism $V_3\ra U_3$ is a covering,
the sheaf axiom  ensures that already $q_\alt^*(\varphi) =0$.
\qed

\medskip
It is easy to see that the class of $(U_\bullet,\shT,\varphi)$ in the group $\RTC^2(\shF)$
is independent of the choice of the
equivalence class of the $\shF$-gerbe $\foG\ra\foC$, the covering $U\ra X$, and the object $T\in\foG_U$.
One thus gets a well-defined map
$$
H^2(X,\shF)\lra \RTC^2(\shF),\quad \foG\longmapsto (U_\bullet,\shT,\varphi).
$$
In the next section, we shall see that this map is bijective.

\section{The comparison map}
\mylabel{Comparison map}

We keep the assumptions of the preceding section,
such that $\foC$ is a site 
having  a final object $X\in\foC$,
and that the Grothendieck topology is given by a pretopology
$\Cov(V)$, $V\in\foC$ of covering families $(U_\lambda\ra V)_{\lambda\in L}$.
Furthermore, we suppose that for each covering family $(U_\lambda\ra V)_{\lambda\in L}$,
the disjoint union $U=\bigcup U_\lambda$ exists, such that each covering family can be refined to 
a covering single. For each abelian sheaf $\shF$, fix an injective
resolution $0\ra\shF\ra\shI^0\ra\ldots$, such that sheaf cohomology becomes
$H^n(X,\shF)=H^n\Gamma(X,\shI^\bullet)$.

The goal of this section is to identify the group $\RTC^n(\shF)$ of equivalence classes of 
rigidified torsor cocycles with the cohomology group $H^n(X,\shF)$, for each degree $n\geq 1$.
The crucial ingredient is the three-term complex $C^{n-1}\stackrel{\Psi}{\ra} C^n\stackrel{\Phi}{\ra} C^{n+1}$ constructed in
Section \ref{Three-term complex}. The main task is to define the \emph{comparison map}
$$
\RTC^n(\shF)\lra \dirlim\Kernel(\Phi)/\Image(\Psi) = H^n(X,\shF),
$$
where the identification on the right comes from Theorem \ref{cohomology three-term}.
We start to define the comparison map 
on  objects  $(U_\bullet,\shT,\varphi)$ from a fixed fiber category $\foR_\shF^n(U_\bullet)$.
Recall that $\foR_\shF^n\ra\foH_{X,n-2}$ is the  fibered Picard category of rigidified  torsor
cocycles.

Consider the  induced $\shI^0|U_{n-1}$-torsor
$\shT^0=(\shI^0|U_{n-1})\wedge^{\shF|U_{n-1}}\shT$, obtained by extending the structure 
sheaf with respect to the inclusion $\shF\subset\shI^0$.
This torsor admits a section $s\in\Gamma(U_{n-1},\shT^0)$, because $\shI^0$ is an injective and hence acyclic
sheaf. The resulting bijection $\shT^0\ra \shI^0|U_{n-1}$ yields an injection $\shT\subset\shI^0|U_{n-1}$,
whose image under the vertical differential $\epsilon d:\shI^0\ra\shI^1$ can be regarded as a section $\alpha_{n-1}\in\Gamma(U_{n-1},\shI^1)$
mapping to zero in $\Gamma(U_{n-1},\shI^2)$. Here $\epsilon=(-1)^{n-1}$ is the sign introduced for the double complex.

Next, consider the section $\varphi\in\Gamma(U_n,p_\alt^*(\shT))$ for the alternating preimage.
Recall that $p_i:U_n\ra U_{n-1}$, $0\leq i\leq n$ denote the face operators in the hypercovering $U_\bullet$.
The inclusion $\shT\subset\shT^0=\shI^0|U_{n-1}$ induces an inclusion of alternating preimages
$$
p_\alt^*(\shT)\subset p_\alt^*(\shT^0) = p_\alt^*(\shI^0|U_{n-1})=\shI^0|U_n,
$$
so our section $\varphi$ becomes an element $\alpha_n\in\Gamma(U_n,\shI^0)$.
In turn, we obtain a cochain   $(\alpha_{n-1},\alpha_n)\in C^n$ in the three-term complex  
$C^{n-1}\stackrel{\Psi}{\ra} C^n\stackrel{\Phi}{\ra} C^{n+1}$, and the comparison map
will be given by the assignment 
\begin{equation}
\label{comparison on objects}
(U_\bullet,\shT,\varphi)\longmapsto (\alpha_{n-1},\alpha_n).
\end{equation}
This cochain is a cocycle: We already remarked above that $d(\alpha_{n-1})=0$. 
By construction, $\alpha_n$ is the section $\varphi$ of the subsheaf $p_\alt^*(\shT)\subset\shI^0|U_n$,
hence $-\epsilon d(\alpha_n)=\delta(\alpha_{n-1})$.
Furthermore, we have $q_\alt^*(\varphi)=0$ as section of the iterated alternating preimage
$q_\alt^*(p_\alt^*(\shT))$, which ensures $\partial(\alpha_n)=0$.
Summing up,  the pair $(\alpha_{n-1},\alpha_n)$ lies in the kernel of the   differential $\Phi$, hence  is a cocycle. 

This attaches to each object $(U_\bullet,\shT,\varphi)\in\foR_\shF^n(U_\bullet)$ an element $(\alpha_{n-1},\alpha_n)\in\Kernel(\Phi)$.
Note that the assignment depends on the choice of sections $s\in\Gamma(U_{n-1},\shT^0)$. Passing to cohomology   
gives a map
\begin{equation}
\label{comparision map objects}
\foR_\shF^n(U_\bullet)\lra \Kernel(\Phi)/\Image(\Psi) \lra H^n(X,\shF),\quad
(U_\bullet,\shT,\varphi)\longmapsto (\alpha_{n-1},\alpha_n).
\end{equation}
 This map does not depend anymore on
the choice of sections: Any other section is of the form $s'=s+\beta_{n-1}$ for
some unique $\beta_{n-1}\in\Gamma(U_{n-1},\shI^0)$. One easy checks that
the resulting cocycle $(\alpha'_{n-1},\alpha_n')$ differs by the coboundary $\Psi(0, \beta_{n-1})$.
The following is also immediate:

\begin{proposition}
\mylabel{properties comparison objects}
The above map \eqref{comparision map objects} sends isomorphic objects   to the same cohomology class,
and turns  wedge products in the Picard category $\foR_\shF^n(U_\bullet)$
into addition of classes. Furthermore, for each $\shF|U_{n-2}$-torsor $\shT_{n-2}$ and each local section $s\in\Gamma(U_{n-1},\shF)$, 
the  rigidified torsor coboundary   $(U_\bullet,\shT,\varphi)$ given by  $\shT=r^*_\alt(\shT_{n-2})$ and $\varphi=p_\alt^*(s)\cdot\varphi_\can$ is send to the zero class.
\end{proposition}

Clearly, our map  is functorial in $U_\bullet$.
Recall that the group $\RTC^n(\shF)$ is the direct limit for the functor
$$
\overline{\foH}_{X,n-2}^\op\lra (\Ab),\quad U_\bullet\longmapsto \pi_0(\foR^n_\shF(U_\bullet))/\pi_0(\foT^{n-1}_\shF(U_\bullet)),
$$
which is defined on the quotient category of hypercoverings of type $r=n-2$.
Passing to direct limits, the maps in \eqref{comparision map objects} give the desired comparison map
\begin{equation}
\label{comparision map}
\RTC^n(\shF)\lra   \dirlim\Kernel(\Phi)/\Image(\Psi) = H^n(X,\shF).
\end{equation}
We now come to our  main result, which gives the desired geometric interpretation of higher cohomology:

\begin{theorem}
\mylabel{comparison bijective}
For each abelian sheaf $\shF$ and each degree $n\geq 1$, the     comparison map  is bijective,
such that we have an identification $\RTC^n(\shF)=H^n(X,\shF)$.
\end{theorem}

\proof
To see surjectivity,  we represent a given cohomology class $[\alpha]\in H^n(X,\shF)$ by a cocycle $(\alpha_{n-1},\alpha_n)$
in the three-term complex $C^{n-1}\stackrel{\Psi}{\ra} C^n\stackrel{\Phi}{\ra} C^{n+1}$, with respect
to some hypercovering $U_\bullet$ of type $r=n-2$.
In particular, the entry $\alpha_{n-1}\in\Gamma(U_{n-1},\shI^1)$ is a local section that vanishes in $\Gamma(U_{n-1},\shI^2)$.
The cartesian square of set-valued sheaves
$$
\xymatrix{
\shT\ar[r]\ar[d]				& 	h_{U_{n-1}}\ar[d]^{\alpha_{n-1}}\\
\shI^0|U_{n-1}\ar[r]_{\epsilon d}		& 	\shI^1|U_{n-1}
}
$$
defines an $\shF|U_{n-1}$-torsor $\shT$.  Here $\epsilon=(-1)^{n-1}$ is the sign introduced for the double complex.
One easily checks that its alternating preimage sits
in the cartesian square
$$
\xymatrix{
p_\alt^*(\shT)\ar[r]\ar[d]	& 	h_{U_n}\ar@{.>}[dl]_{-\alpha_n}\ar[d]^{\partial(\alpha_{n-1})}\ar@{.>}@/_1.5em/[l]_{\varphi}\\
\shI^0|U_n\ar[r]_{\epsilon d}		& 	\shI^1|U_n.
}
$$
The   dotted diagonal arrow arises from the   the entry $\alpha_n\in\Gamma(U_n,\shI^0)$, 
in light of the cocycle condition $\partial(\alpha_{n-1})=-\epsilon d(\alpha_{n})$. 
This dotted diagonal arrow corresponds to a section $\varphi\in\Gamma(U_n,p_\alt^*(\shT))$. In turn, we obtain a triple $(U_\bullet,\shT,\varphi)$.
The condition $\partial(\alpha_n)=0$ ensures that $q_\alt^*(\varphi)=0$, hence $(U_\bullet,\shT,\varphi)$ is a rigidified torsor
cocycle. 

We now check that  $(U_\bullet,\shT,\varphi)\mapsto (\alpha_n,\alpha_{n-1})$ under the comparison map, as described on objects in   
\eqref{comparison on objects}.
Indeed, the inclusion $\shT\subset\shI^0|U_{n-1}$ yields a canonical isomorphism
$$
\shT^0=(\shI^0|U_{n-1})\wedge^{\shF|U_{n-1}}\shT\lra \shI^0|U_{n-1},
$$
and we may take for the section $s\in\Gamma(U_{n-1},\shT^0)$   the zero-section of $\shI^0|U_{n-1}$. Then the image of $\shT\subset\shI^0|U_{n-1}$
under the differential $\epsilon d$ is the entry $\alpha_{n-1}\in\Gamma(U_{n-1},\shI^1)$. Furthermore, the section $\varphi$
corresponds to $-\alpha_n\in\Gamma(U_n,\shI^0)$. Thus the cocycle $(\alpha_{n-1},\alpha_n)$ lies in the image.

It remains to check that the comparison map is injective.
Let $(U_\bullet,\shT,\varphi)$ be a rigidified torsor cocycle, and suppose that the comparison map 
\eqref{comparison on objects} sends it to a cocycle $(\alpha_{n-1},\alpha_n)$
in the three-term complex $C^{n-1}\stackrel{\Psi}{\ra} C^n\stackrel{\Phi}{\ra} C^{n+1}$ that vanishes in $H^n(X,\shF)$.
Refining the hypercovering $U_\bullet$, we may assume that $(\alpha_{n-1},\alpha_n)$ is already the coboundary of some $(\beta_{n-2},\beta_{n-1})$,
in other words
$$
\alpha_{n-1}= \epsilon d(\beta_{n-1}) + \sum_{k=0}^{n-1} (-1)^k r_k^*(\beta_{n-2}) \quadand
\alpha_n =  \sum_{i=0}^n (-1)^i p_i^*(\beta_{n-1}).
$$
The local section $\beta_{n-2}\in\Gamma(U_{n-2},\shI^0/\shF)$ yields an $\shF|U_{n-2}$-torsor $\shT_{U_{n-2}}$,
by taking the preimage sheaf of this local section.  Adding the rigidified torsor coboundary $(U_\bullet, r^*_\alt(\shT^{-1}), \varphi_\can)$,
we reduce to the situation $\beta_{n-2}=0$, hence  $\alpha_{n-1}= \epsilon d(\beta_{n-1})$.
By definition of the comparison map, $\shT\subset\shI^0|U_{n-1}$ is the preimage sheaf of $\alpha_{n-1}\in\Gamma(U_{n-1},\shI^1)$,
so $\epsilon\beta_{n-1}\in\Gamma(U_{n-1},\shI^0)$ defines a section for $\shT$.
In turn, we may assume that $\shT=\shF|U_{n-1}$ and $\beta_{n-1}=0$. Consequently $\alpha_n=0$,
which means that the section $\varphi\in\Gamma(U_{n+1},q_\alt^*(\shT))$ coincides with the zero section.
Summing up, our rigidified torsor cocycle is a rigidified torsor coboundary.
\qed

\medskip
Let us record the following  two consequences:  

\begin{corollary}
\mylabel{same cohomology class}
Two rigidified   $\shF$-torsor  $n$-cocycles $A'=(U'_\bullet,\shT',\varphi')$ and $A''=(U''_\bullet,\shT'',\varphi'')$
have the same cohomology class in $H^n(X,\shF)$ if and only if there is a common refinement 
$U'_\bullet \leftarrow U_\bullet\ra U''_\bullet$  
and some rigidified torsor coboundary $B=(U_\bullet, r_\alt^*(\shT_{n-2}),p_\alt^*(s)\cdot\varphi_\can)$
such that  $A'|U_\bullet \simeq (A''|U_\bullet) \wedge B$.
\end{corollary}

\begin{corollary}
\mylabel{trivial cohomology class}
A rigidified   $\shF$-torsor  $n$-cocycles $A'=(U'_\bullet,\shT',\varphi')$ has trivial
cohomology class in $H^n(X,\shF)$ if and only if there is a refinement 
$U_\bullet\ra U'_\bullet$  
and some rigidified torsor coboundary $B=(U_\bullet, r_\alt^*(\shT_{n-2}),p^*_\alt(s)\cdot\varphi_\can)$
such that  $A'|U_\bullet \simeq   B$.
\end{corollary}

\section{Bundle gerbes and Dixmier--Douady classes}
\mylabel{Bundle gerbes}

Let us now connect our theory of rigidified torsor cocycles with 
Murray's notion  of bundle gerbes from \cite{Murray 1996}, Section 3, see also  \cite{Murray; Stevenson 2011}, Section 2. 
Let $M$ be a  differential manifold,   $\pi:Y\ra M$ be a fibration, 
and   $P\ra Y^{[2]}$ be a $\CC^\times$-principal bundle.
In this context, $Y$ is also a differentiable manifold  and   fibrations
denote  differentiable maps that are surjective on points and tangent vectors.
The total space $Y$ is allowed to be infinite-dimensional,   and the fibration $Y\ra M$ is assumed to admit local sections.
Conforming with  the notation in loc.\ cit., we here write $Y^{[2]}=Y\times_MY$ for the fiber product.
Now define 
$$
Y^{[2]}\circ Y^{[2]}\subset Y^{[2]}\times Y^{[2]}
$$
as the set of all pairs
of the form $((a,b),(b,c))$, with $a,b,c\in Y$ all mapping to the same point in $M$.  
Write 
$\pi_i:Y^{[2]}\circ Y^{[2]}\ra Y^{[2]}$ with  $i=1,2$ for the two projections.

A \emph{bundle gerbe} consists of the choice of a fibration $\pi:Y\ra M$
and a principal $\CC^\times$-bundle $P\ra Y^{[2]}$, together with a
map of $\CC^\times$-bundles 
\begin{equation}
\label{product}
\mu:\pi_1^{-1}(P)\otimes\pi_2^{-1}(P)\lra P
\end{equation}
covering the   map 
\begin{equation}
\label{triple}
\pi_3:Y^{[2]}\circ Y^{[2]}\lra Y^{[2]},\quad  ((a,b),(b,c))\longmapsto (a,c).
\end{equation}
The product in \eqref{product} is assumed to be \emph{associative} whenever triple products in \eqref{triple} are defined. 
Bundle gerbes arising from $P=\pi_1^{-1}(Q^*)\otimes\pi_2^{-1}(Q)$, where $Q\ra Y$ is a principal $\CC^\times$-bundle, are called \emph{trivial}.

Let us translate this into the set-up and notation of the present paper.
The tensor product in \eqref{product} is    vector bundle notation for the contracted product
of $\CC^\times$-bundles. Furthermore, the morphism can be regarded as a map
$\mu:\pi_1^*(P)\otimes\pi_2^*(P)\ra \pi_3^*(P)$ of principal $\CC^\times$-bundles over $Y^{[2]}\circ Y^{[2]}$.
Using semi-simplicial notation, we write  $Y_d=Y^{[d+1]}$. Clearly, the canonical map
$$
Y^{[2]}\circ Y^{[2]}\lra Y^{[3]}= Y_2,\quad ((a,b),(b,c))\longmapsto (a,b,c)
$$
into the threefold fiber product is a homeomorphism. With respect to this identification, we have 
$$
\pi_1=p_2\quadand \pi_2=p_0\quadand \pi_3=p_1
$$
as face operators $Y_2\ra Y_1$, in simplicial notation. So  $\mu$ in  \eqref{product}  may be regarded a section
into the alternating preimage $p_\alt^*(P^{-1})=p_0^*(P^{-1})\otimes p_1^*(P)\otimes p_2^*(P^{-1})$ 
for the inverse principal bundle $P^{-1}$. On the other hand, we may regard
$\mu$ as the collection of fiber-wise maps of principal $\CC^\times$-sets 
$$
\mu_{abc}:P_{(a,b)}\otimes P_{(b,c)}\lra P_{(a,c)},
$$
where $a,b,c\in Y$ map to the same point in $M$. Write $\id_{ab}:P_{(a,b)}\ra P_{(a,b)}$ for the identity map. 
The associativity condition for bundle gerbes becomes
\begin{equation}
\label{associative}
\mu_{acd}\circ (\mu_{abc}\otimes\id_{cd})= \mu_{abd} \circ ( \id_{ab}\otimes\mu_{bcd})
\end{equation}
as maps $P_{(a,b)}\otimes P_{(b,c)}\otimes P_{(c,d)}\ra P_{(a,d)}$ of principal $\CC^\times$-sets, 
for all $a,b,c,d\in Y$ mapping to the same point in $M$. For the following observation, recall 
that the $p_i:Y_2\ra Y_1$ and $q_j:Y_3\ra Y_2$ denote   face operators.

\begin{proposition}
The associativity condition \eqref{associative} for bundle gerbes holds if and only if $q_\alt^*(\mu)=1_{Y_3}$
with respect to the identification $q_\alt^* (p_\alt^*(P^{-1}))=\CC^\times\times Y_3$ of iterated alternating preimages.
\end{proposition}

\proof
Over each point $(a,b,c,d)\in Y^{[4]}=Y_3$, the fiber of the iterated alternating preimage $q_\alt^*(p_\alt^*(P^{-1}))$
is the tensor product of the following twelve principal $\CC^\times$-sets:
$$
\begin{array}{ccc}
P_{(c,d)}^{-1}	&	P_{(b,d)}	& 	P_{(b,c)}^{-1}\\
P_{(c,d)}	&	P_{(a,d)}^{-1}	& 	P_{(a,c)}\\
P_{(b,d)}^{-1}	&	P_{(a,d)}	& 	P_{(a,b)}^{-1}\\
P_{(b,c)}	&	P_{(a,c)}^{-1}	& 	P_{(a,b)}\\
\end{array}
$$
Now choose for each of the six occurring points   $(c,d),\ldots,(a,b)\in Y^{[2]}$ some elements $s_{cd}\in P_{(c,d)}, \ldots, s_{ab}\in P_{(a,b)}$.
Regarding $\mu$ as a pairing, the four equations
\begin{gather*}
\mu_{bcd}(s_{bc}\otimes s_{cd}) = f_a s_{bd},\\
\mu_{acd}(s_{ac}\otimes s_{cd}) = f_b s_{ad},\\
\mu_{abd}(s_{ab}\otimes s_{bd}) = f_c s_{ad},\\
\mu_{abc}(s_{ab}\otimes s_{bc}) = f_d s_{ac} 
\end{gather*}
define scalars $f_a,\ldots,f_d\in \CC^\times$.
The condition $q_\alt^*(\mu)=1_{Y_3}$ on the alternating preimage  translates into $f_af_c=f_bf_d$. Note that here
we use multiplicative rather than additive notation.
Applying the two sides of \eqref{associative} to the element 
$$
s_{ab}\otimes s_{bc}\otimes s_{cd}\in P_{(a,b)}\otimes P_{(b,c)}\otimes P_{(c,d)},
$$
we see that the associativity conditions is   equivalent to $f_af_c=f_bf_d$ as well.
\qed

\medskip
It is now straightforward to verify that Murray's bundle gerbes correspond to 
our rigidified torsor cocycles for $n=2$:
Let $\shF=\shC^\times_M$ be the abelian sheaf of invertible complex-valued functions.
Each principal $\CC^\times$-bundle $P\ra M$ yields the  $\shF$-torsor  $\shT$ of local sections,
and each $\shF$-torsor $\shT$ yields the principal $\CC^\times$-bundle defined as the
relative spectrum of the $\shC_M$-algebra $\bigoplus_{d\in\ZZ}\shL^{\otimes -d}$, 
where $\shL=\shC_M\wedge^{\shC^\times_M}\shT$ is the invertible sheaf attached to the $\shC_M^\times$-torsor.
In this way, one gets an equivalence of categories between principal $\CC^\times$-bundles 
and $\shC_M^\times$-torsors.

In turn, a bundle gerbe, which consists of a fibration $P\ra M$,
a principal bundle $P\ra Y^{[2]}$ and a pairing $\mu:\pi_1^{-1}(P)\otimes\pi_2^{-1}(P)\ra P$
satisfying the associativity condition,
corresponds to a covering $U=Y$ with respect to some suitable site $\foC$, 
an $\shF|U_1$-torsor $\shT$ and a section $\varphi\in\Gamma(U_2,p_\alt^*(\shT))$
satisfying $q_\alt^*(\varphi)=1$, in multiplicative notation.

The Dixmier--Douady class arises as follows: The exponential sequence of abelian sheaves
$0\ra 2\pi i\ZZ_M\ra\shC_M\ra\shC_M^\times\ra 1$ induces   long exact sequences
$$
H^n(M,\shC_M)\lra H^n(X,\shC_M^\times)\lra H^n(M,2\pi i\ZZ)\lra H^{n+1}(M,\shC_M).
$$
The outer terms vanish for $n\geq 1$, because the  sheaf $\shC_M$ is soft  whence acyclic.
In turn, we get the identifications
$$
\RTC^n(\shC_M^\times) = H^n(M,\shC_M^\times) = H^{n+1}(M,2\pi i\ZZ).
$$
Here one may interpret the right hand side both as sheaf cohomology and singular cohomology.
The integral cohomology class attached to a bundle gerbe or a rigidified torsor cocycle
in degree $n=2$ is called the \emph{Dixmier--Douady class}. We see that
it is defined for all degree's $n\geq 1$.

Note also the theory   works well if the differentiable manifold $M$ is merely a  topological space that is paracompact and 
locally contractible, because then sheaf cohomology with locally constant $\ZZ$-coefficients
coincides with singular cohomology, as explained in \cite{Bredon 1967}, Chapter III.


\end{document}